\documentclass{article}
\usepackage[utf8]{inputenc}

\pdfoutput=1
\usepackage{mathtools}
\usepackage{amsfonts}
\usepackage{amssymb}
\usepackage{centernot}
\usepackage{graphicx}
\usepackage{amsthm}
\usepackage{enumerate}

\DeclareMathOperator{\coker}{coker }
\usepackage{tikz}
\usepackage{tikz-cd}
\usetikzlibrary{decorations.markings}
\usetikzlibrary{arrows}
\usetikzlibrary{calc}
\providecommand{\U}[1]{\protect\rule{.1in}{.1in}}
\newtheorem{theorem}{Theorem}

\newtheorem{corollary}[theorem]{Corollary}

\newtheorem{definition}[theorem]{Definition}

\newtheorem{lemma}[theorem]{Lemma}

\newtheorem{proposition}[theorem]{Proposition}

\theoremstyle{remark}

\newcommand{\MQ}{\textup{MQ}}
\newcommand{\IMQ}{\textup{IMQ}}

\newcommand{\Aut}{\textup{Aut}}
\newcommand{\Dis}{\textup{Dis}}
\newcommand{\Fix}{\textup{Fix}}

\newcommand{\Qr}{Q_A^{\textup{red}}}
\newcommand{\Mr}{M_A^{\textup{red}}}
\newcommand{\Me}{M_A^{\textup{enr}}}
\newcommand{\fr}{\phi_L^{\textup{red}}}
\newcommand{\ann}{\textup{ann}}

\newcommand{\inv}{\textup{inv}}

\allowdisplaybreaks

\begin{document}

\title{Multivariate Alexander quandles, V. Constructing the medial quandle of a link}
\author{Lorenzo Traldi\\Lafayette College\\Easton, PA 18042, USA\\traldil@lafayette.edu
}
\date{ }
\maketitle

\begin{abstract}
We explain how the medial quandle of a classical or virtual link can be built from the peripheral structure of the reduced Alexander module.


\end{abstract}

\section{Introduction}

This is the fifth paper in a series involving Alexander modules and quandles, but the paper is written so that familiarity with the earlier papers in the series is not required. The purpose of the paper is to show that the medial quandle of a link can be constructed from the reduced Alexander module. The construction requires the module's peripheral elements, which were introduced recently \cite{peri}. In the case of a knot, the peripheral elements do not contribute anything to the construction, and we deduce an extension to virtual knots of Joyce's result \cite[Theorem 17.3]{J} that the medial quandle of a classical knot is the standard Alexander quandle on the Alexander invariant of the knot. 

To keep the introduction as short as possible, we leave the statements of many definitions and previous results for later in the paper. We should say, though, that we use the terms \emph{knots} and \emph{links} for oriented virtual knots and links. 

Let $\Lambda=\mathbb Z[t^{\pm 1}]$ be the ring of Laurent polynomials in the variable $t$, with integer coefficients. Here is the central algebraic construction of the paper.

\begin{definition}
\label{qquandle}
Suppose $M$ is a $\Lambda$-module with a submodule $N$, $I$ is a nonempty set, $m_i \in M\thickspace \allowbreak \forall i \in I$, and $m_i-m_j \in N \thickspace \allowbreak \forall i,j \in I$. Suppose also that for each $i \in I$, $X_i$ is a submodule of $N$ such that $(1-t) \cdot X_i = 0$. Then there is an associated medial quandle $Q(N,(m_i), \allowbreak (X_i))$, defined as follows.
\begin{enumerate}
    \item For each $i \in I$, let $Q_i=N/X_i$. If $i \neq j \in I$ then the sets $Q_i$ and $Q_j$ are understood to be disjoint. Define the set $Q(N,(m_i), \allowbreak (X_i))$ by
    \[
    Q(N,(m_i), \allowbreak (X_i))=\bigcup_{i \in I} Q_i .
    \]
    \item Define an operation $\triangleright$ on $Q(N,(m_i), \allowbreak (X_i))$ as follows. If $i,j \in I$, $x\in Q_i$ and $y \in Q_j$, then
    \[
    x \triangleright y = m_j-m_i+tx+(1-t)y+X_i \in Q_i.
    \]
\end{enumerate}
\end{definition}

Here are three comments on Definition \ref{qquandle}. First: the formula for $\triangleright$ in part 2 of Definition \ref{qquandle} is well defined, because the fact that $(1-t) \cdot X_j=0$ guarantees there is no ambiguity in the value of the coset $(1-t)y+X_i$. Second: Definition \ref{qquandle} allows $M$ to strictly contain $N$, but this is a merely stylistic choice. If we would like to have $M=N$, we can always replace the elements $m_i \in M$ with elements $m'_i \in N$, by choosing a fixed index $i_0 \in I$ and defining $m'_i = m_i-m_{i_0} \thickspace \allowbreak \forall i \in I$. It is easy to see that then $Q(N,(m_i), \allowbreak (X_i)) \cong Q(N, (m'_i), (X_i))$. Third: If $\mu$ is a positive integer and $I=\{1, \dots, \mu\}$, we use the notation $Q(N,(m_i), \allowbreak (X_i))=Q(N,m_1, \dots, m_\mu, X_1, \dots, X_\mu)$.

In Sec.\ \ref{def1} we show that Definition \ref{qquandle} provides a general description of medial quandles. To be precise:

\begin{theorem}
\label{main1}
If $Q$ is a medial quandle, then there are $M,N,I, (m_i)$ and $(X_i)$ such that $Q \cong Q(N,(m_i),(X_i))$.
\end{theorem}

We should mention that Definition \ref{qquandle} and Theorem \ref{main1} are closely related to the structure theory of medial quandles due to Jedli\v{c}ka \emph{et al}.\ \cite{JPSZ1, JPSZ2}. In particular,  Theorem \ref{main1} extends and refines \cite[Lemma 4.3]{JPSZ2}.

The reduced Alexander module $\Mr(L)$ is a well-known invariant of a link $L$. For a classical link, it corresponds to the first relative homology group of the total linking number cover, with respect to its fiber. The module is conveniently described by a presentation with generators and relations corresponding to the arcs and crossings of a link diagram (respectively). There is a $\Lambda$-linear map $\fr:\Mr(L) \to \Lambda$,  under which $1$ is the image of every module generator corresponding to an arc in a diagram. The kernel of $\fr$ is the \emph{reduced Alexander invariant} of $L$.

In recent work \cite{peri}, we introduced peripheral elements in the reduced Alexander module of a link. If $L=K_1 \cup \dots \cup K_ \mu$ is a $\mu$-component link then there are subsets $M_1(L), \dots, \allowbreak M_ \mu(L) \subset \Mr(L)$; the elements of $M_i(L)$ are the \emph{meridians of $K_i$ in $\Mr(L)$}. In a module presentation of $\Mr(L)$ corresponding to a diagram of $L$, every generator corresponding to an arc of $K_i$ is a meridian associated with $K_i$; also, every meridian is mapped to $1$ by $\fr$. For each $i \in \{1, \dots, \mu \}$ there is a single longitude $\chi_i(L)$ corresponding to $K_i$. These longitudes have several distinctive properties, including $\chi_i(L) \in \ker \fr$ and $(1-t) \cdot \chi_i(L)=0$. The list
\[
(\Mr(L),M_1(L), \dots, M_\mu(L),\chi_1(L), \dots , \chi_\mu(L))
\]
is the \emph{enhanced reduced Alexander module} of $L$; it is denoted $\Me(L)$. The enhanced reduced Alexander module is a link invariant, in the sense that if $L$ and $L'$ are equivalent links then there is a $\Lambda$-linear isomorphism $\Mr(L) \to \Mr(L')$ such that for each $i \in \{1, \dots, \mu \}$, the image of $M_i(L)$ is $M_i(L')$ and the image of  $\chi_i(L)$ is $\chi_i(L')$.

Here is the main theorem of this paper. 

\begin{theorem}
\label{main}
Let $L$ be a link. For each $i \in \{1, \dots, \mu \}$, let $m_i \in M_i(L)$ be any meridian of $K_i$ in $\Mr(L)$, and let $X_i$ be the cyclic submodule of $\Mr(L)$ generated by $\chi_i(L)$. Then the quandle $Q(\ker \fr, m_1, \dots, m_\mu, X_1, \dots, X_\mu)$ of Definition \ref{qquandle} is isomorphic to the medial quandle of $L$.
\end{theorem}

We should mention that the medial quandle of $L$ was denoted $\textup{AbQ}(L)$ in Joyce's seminal work \cite{J}. We use the notation $\MQ(L)$ instead.

In the special case when $L$ is a knot (i.e.\ $\mu=1$), the longitude $\chi_1(L)$ is $0$ \cite{peri}. Theorem \ref{main} then states that $\MQ(L)$ is the quandle on $\ker \fr$ given by the operation $x \triangleright y = tx + (1-t)y$. This is the standard way to define a quandle structure on a $\Lambda$-module; in the literature such a quandle is called an affine quandle or an Alexander quandle. We deduce the following extension of a classical result of Joyce \cite[Theorem 17.3]{J} to virtual knots. 

\begin{corollary} 
\label{knotcor}
Suppose $L$ is a knot. Then the medial quandle of $L$ is isomorphic to the standard Alexander quandle on the $\Lambda$-module $\ker \fr$. Moreover, $\MQ(L)$ and $\ker \fr$ are equivalent as invariants of $L$, i.e.\ each invariant determines the other.
\end{corollary}

In general, Theorem \ref{main} implies that the medial quandle of a link is always determined by the enhanced reduced Alexander module $\Me(L)$. When $\mu>1$, $\Me(L)$ is a strictly stronger link invariant than $\MQ(L)$. There are two simple reasons for this. The first reason is that $\Me(L)$ involves the indexing of the components of $L= \allowbreak K_1 \cup \dots \cup K_ \mu$, and the medial quandle does not. Therefore if $L'$ is obtained from $L$ by permuting the indices of $K_1, \dots, K_\mu$ then $\Me(L)$ might distinguish $L$ from $L'$, but the medial quandle cannot. The second reason is that $\Me(L)$ involves the individual longitudes $\chi_1(L), \dots, \chi_ \mu(L)$, while $Q(\ker \fr, m_1, \dots, m_\mu, \allowbreak X_1, \dots, X_\mu)$ involves only the cyclic submodules $X_1, \dots, X_\mu$ generated by $\chi_1(L), \dots, \chi_ \mu(L)$. Therefore if $\Me(L)$ and $\Me(L')$ differ only because their longitudes are negatives of each other, then $\Me(L)$ might distinguish $L$ from $L'$, but the medial quandle cannot. 

Here is an outline of the paper. In Sec.\ \ref{mq}, we briefly summarize some aspects of the general theory of medial quandles. In Sec.\ \ref{def1}, we discuss the special properties of the quandles given by Definition \ref{qquandle}, and prove Theorem \ref{main1}. Notation for link diagrams and link invariants is established in Sec.\ \ref{linknot}, and Theorem \ref{main} is proven in Sec.\ \ref{mainproof}. In Sec.\ \ref{semireg} we observe that the quandle $\Qr(L)$ of \cite{mvaq4} is the maximal semiregular image of $\MQ(L)$. In Sec.\ \ref{invmed} we apply our discussion to the involutory medial quandle $\IMQ(L)$, sharpening some results from \cite{mvaq2}. Several examples are presented in Sec.\ \ref{exsec}.

Before proceeding we should thank an anonymous reader for helpful comments and corrections.

\section{Medial Quandles}
\label{mq}

In this section we present some properties of medial quandles, without much explanation or proof. For a more detailed and thorough account of the theory we refer to Jedli\v{c}ka \emph{et al}.\ \cite{JPSZ1, JPSZ2}. Much of what we summarize here is also detailed in the preceding paper in this series \cite{mvaq4}. Before jumping in, we should mention that in his seminal work \cite{J}, Joyce used the term ``abelian'' rather than ``medial.''

\begin{definition}
\label{mqdef}
A \emph{medial quandle} is given by a binary operation $\triangleright$ on a set $Q$. The following properties must hold.
\begin{enumerate} 
\item $x\triangleright x=x \thickspace \allowbreak \forall x \in Q$.
\item For each $y \in Q$, a permutation $\beta_y$ of $Q$ is defined by $\beta_y(x)=x \triangleright y$.
\item $(w\triangleright x) \triangleright (y \triangleright z)=(w\triangleright y) \triangleright (x\triangleright z) \thickspace \allowbreak \forall w,x,y,z \in Q$.
\end{enumerate}
\end{definition}

All the quandles we consider in this paper are medial, but it is worth mentioning that a general quandle is required to satisfy the first two properties of Definition \ref{mqdef}, along with the special case of the third property in which $y=z$. 

If $Q$ is a medial quandle then the maps $\beta_y:Q \to Q$ are the \emph{translations} or \emph{inner automorphisms} of $Q$. We often use the notation $\beta_y^{-1}(x)=x \triangleright^{-1} y$. The translations are automorphisms of $Q$, and the subgroup of $\Aut(Q)$ they generate is denoted $\beta(Q)$. The composition $\beta_y \beta ^{-1} _z$ of a translation with the inverse of a translation is an \emph{elementary displacement} of $Q$, and the subgroup of $\beta(Q)$ generated by the elementary displacements is the \emph{displacement group} $\Dis(Q)$. Note that the inverse of an elementary displacement is also an elementary displacement, so every displacement is a composition of elementary displacements. In fact it turns out that $\Dis(Q)$ includes all products $\prod \beta_{y_i}^{m_i}$ with $\sum m_i=0$. As discussed in \cite[Sec.\ 5]{mvaq4}, the displacement group is a normal abelian subgroup of $\beta(Q)$, and it may be given the structure of a $\Lambda$-module by choosing any fixed element $q^* \in Q$, and defining $t\cdot d = \beta_{q^*} d \beta_{q^*}^{-1} \thickspace \allowbreak \forall d \in \Dis(Q)$. Changing the choice of $q^*$ does not change the $\Lambda$-module structure of $\Dis(Q)$. The reason is simple: if $q^*,q^{**} \in Q$ then $\beta_{q^*}\beta^{-1}_{q^{**}} \in \Dis(Q)$, so commutativity of $\Dis(Q)$ implies the following.
\begin{align*}
    \beta_{q^*}d\beta_{q^*}^{-1}
   & =\beta_{q^*}d\beta_{q^*}^{-1}(\beta_{q^*}\beta^{-1}_{q^{**}})(\beta_{q^*}\beta^{-1}_{q^{**}})^{-1}\\
   &=(\beta_{q^*}\beta^{-1}_{q^{**}})^{-1}\beta_{q^*}d\beta_{q^*}^{-1}(\beta_{q^*}\beta^{-1}_{q^{**}})\\
   &=\beta_{q^{**}}\beta_{q^*}^{-1}\beta_{q^*}d\beta_{q^*}^{-1}\beta_{q^*}\beta^{-1}_{q^{**}}=\beta_{q^{**}}d\beta^{-1}_{q^{**}}
\end{align*}

Here is an easy consequence of the fact that $t \cdot d = \beta_{q^*}d\beta_{q^*}^{-1}$ does not depend on the choice of $q^*$.

\begin{lemma}
\label{fix}
If $d \in \Dis(Q)$ has a fixed point, then $d = t \cdot d$.
\end{lemma}
\begin{proof}
If $d(q^*)=q^*$ then for every $x \in Q$, we have 
\[
d\beta_{q^*}(x) = d(x \triangleright q^*) = d(x) \triangleright d(q^*) = d(x) \triangleright q^* = \beta_{q^*}d(x).
\]
Therefore $d\beta_{q^*}=\beta_{q^*}d$, so $d=\beta_{q^*}d\beta_{q^*}^{-1}$.
\end{proof}

If $Q_1$ and $Q_2$ are quandles then a function $f:Q_1 \to Q_2$ is a \emph{quandle map} or \emph{quandle homomorphism} if the equality $f(x \triangleright y)=f(x) \triangleright f(y)$ is always satisfied. (That is, $f  \beta_y=\beta_{f(y)}  f \thickspace \allowbreak \forall y \in Q$.) This equality implies
\[
f(x \triangleright^{-1}y) \triangleright f(y) = f((x \triangleright^{-1}y) \triangleright y) = f(x) \text{,}
\]
so the equality $f(x \triangleright^{-1} y)=f(x) \triangleright^{-1} f(y)$ is also always satisfied. (That is, $f  \beta^{-1}_y=\beta^{-1}_{f(y)}  f \thickspace \allowbreak \forall y \in Q$.) The displacement groups of medial quandles are functorial with respect to surjective homomorphisms, in the sense that a surjective quandle map $f:Q_1 \to Q_2$ induces a surjective homomorphism $\Dis(f):\Dis(Q_1) \to \Dis(Q_2)$ of abelian groups, defined in the natural way: If $y_1,\dots,y_n \in Q_1$, $m_1, \dots, m_n \in \{ \pm 1 \}$ and $\sum m_i = 0$, then
    \[
    \Dis(f) \left(\prod \beta_{y_i}^{m_i} \right) =  \prod \beta_{f(y_i)}^{m_i}.
    \]
Note that if $\prod \beta_{y_i}^{m_i} = d \in \Dis(Q_1)$ then as $f$ is a quandle map, 
\begin{align*}
\Dis(f)(d) \circ f &=  \beta_{f(y_1)}^{m_1} \cdots  \beta_{f(y_{n-1})}^{m_{n-1}}  \beta_{f(y_n)}^{m_n} f
\\
&= \beta_{f(y_1)}^{m_1}  \cdots  \beta_{f(y_{n-1})}^{m_{n-1}}  f  \beta_{y_n}^{m_n} 
\\
&= \cdots = f  \beta_{y_1}^{m_1}  \cdots \beta_{y_{n-1}}^{m_{n-1}}  \beta_{y_n}^{m_n} = f \circ d.
\end{align*}

If $q \in Q$ then the \emph{orbit} of $q$ in $Q$ is the smallest subset of $Q$ that contains $q$ and is closed under the action of $\beta(Q)$. We denote the orbit $Q_q$. It turns out that $Q_q = \{d(q) \mid d \in \Dis(Q) \}$. In fact, if $\Fix(q)=\{d \in \Dis(Q) \mid d(q)=q\}$ then $\Fix(q)$ is a $\Lambda$-submodule of $\Dis(Q)$, and $Q_q$ is isomorphic to the standard Alexander quandle on the quotient module $\Dis(Q)/\Fix(q)$. That is, the quandle operation is given by $x \triangleright y = tx +(1-t)y$. An isomorphism can be defined in the natural way: if $d \in \Dis(Q)$ then $d(q) \in Q_q$ corresponds to $d+\Fix(q) \in \Dis(Q)/\Fix(q)$. In particular, a medial quandle $Q$ with only one orbit is isomorphic to the standard Alexander quandle on the $\Lambda$-module $\Dis(Q)$.

As a consequence of these properties, we have the following characterization of isomorphisms of medial quandles.

\begin{proposition}
\label{miso}
Let $Q$ and $Q'$ be medial quandles. Then a quandle map $f:Q \to Q'$ is an isomorphism if and only if it has these three properties.
\begin{enumerate}
    \item The map $f$ is surjective.
    \item Whenever $x,y \in Q$ have different orbits, $f(x),f(y) \in Q'$ also have different orbits.
    \item Whenever $x \in Q$, $d \in \Dis(Q)$ and $\Dis(f)(d)(f(x))=f(x)$, it is also true that $d(x)=x$.
\end{enumerate}
\end{proposition}

\begin{proof}
If $f$ is an isomorphism then the three listed properties certainly hold. For the converse, we need to prove that the three properties imply that $f$ is injective.

Suppose the three properties hold, and $f(x)=f(y)$. The second property implies that $x$ and $y$ have the same orbit. It follows that there is a displacement $d \in \Dis(Q)$ with $d(x)=y$. Then $\Dis(f)(d)(f(x))=(f \circ d)(x)=f(y)=f(x)$, so according to the third listed property, $y=d(x)=x$. 
\end{proof}

One more property of medial quandles will be useful.

\begin{lemma}
\label{disqlem}
Suppose $Q$ is a medial quandle and $S\subseteq Q$ is a generating set. Then the $\Lambda$-module $\Dis(Q)$ is generated by the elements $\beta_s \beta^{-1}_{s'}$ with $s,s' \in S.$
\end{lemma}
\begin{proof}Note first that if $y,z \in Q$ then for every $x \in Q$,
\begin{align*}
\beta_z \beta_{y \triangleright^{-1} z}(x) &=(x \triangleright (y \triangleright^{-1} z)) \triangleright z =(x \triangleright (y \triangleright^{-1} z)) \triangleright (z \triangleright z) \\
&=(x \triangleright z) \triangleright ((y \triangleright^{-1} z) \triangleright z) = (x \triangleright z) \triangleright y = \beta_y \beta_z (x).
\end{align*}
Therefore $\beta_z\beta_{y \triangleright^{-1} z}=\beta_y \beta_z$, and hence $\beta_{y \triangleright^{-1} z} =\beta^{-1}_z \beta_y \beta_z $. If we replace $y \triangleright^{-1} z$ with $w$ in this equality we obtain $\beta_w =\beta_z^{-1} \beta_{ w \triangleright z} \beta_z$, which implies $\beta_{ w \triangleright z}=\beta_z \beta_w \beta_z^{-1}$. Using the two equalities $\beta_{y \triangleright^{-1} z} =\beta^{-1}_z \beta_y \beta_z $ and $\beta_{ w \triangleright z}=\beta_z \beta_w \beta_z^{-1}$, we can use induction to prove that if $q$ and $r$ can be obtained from elements of $S$ using $\triangleright$ and $\triangleright^{-1}$, then $\beta_q \beta_r^{-1}$ can be written as a product of terms $\beta_s^{\pm 1}$ with $s \in S$, such that the sum of the exponents appearing in the product is $0$. 

To complete the proof of the lemma, then, it suffices to show that if $n_1,\dots, \allowbreak n_{k} \in \{\pm 1\}$, $s_1, \dots, s_{k} \in S$ and $\sum n_i = 0$ then $\prod _i \beta_{s_i}^{n_i} \in M_S$, where $M_S$ is the $\Lambda$-submodule of $\Dis(Q)$ generated by $\{ \beta_s \beta^{-1}_{s'} \mid s,s' \in S\}$. We consider three cases. 

If $n_1 \neq n_2$, then 
\[
\prod _{i=1}^{k} \beta_{s_i}^{n_i} = (\beta_{s_1}^{n_1}\beta_{s_2}^{n_2}) \left(\prod _{i=3}^{k} \beta_{s_i}^{n_i} \right) \text{,}
\]
the sum of $\beta_{s_1}^{n_1}\beta_{s_2}^{n_2}$ and $\prod _{i \geq 3} \beta_{s_i}^{n_i}$ in the $\Lambda$-module $\Dis(Q)$. If $n_1=1$, the first summand is an element of $M_S$ by definition. Suppose $n_1=-1$. Recall that $t \cdot d = \beta_{q^*} d \beta^{-1}_{q^*}$ for each $d \in \Dis(Q)$ and any $q^* \in Q$. It follows that $t^{-1} \cdot d = \beta_{q^*}^{-1} d \beta_{q^*}$  for each $d \in \Dis(Q)$ and any $q^* \in Q$, so the first summand is
\[
\beta_{s_1}^{n_1}\beta_{s_2}^{n_2}= \beta^{-1}_{s_1} \beta_{s_2} = \beta^{-1}_{s_1} (\beta_{s_2} \beta^{-1}_{s_1}) \beta_{s_1} = t^{-1} \cdot (\beta_{s_2} \beta^{-1}_{s_1})\text{,}
\]
which is an element of $M_S$. If $k>2$, induction tells us the second summand is also an element of $M_S$. 

If $j>1$, $n_i=-1 \thickspace \allowbreak \forall i \in \{1, \dots, j\}$ and $n_{j+1}=1$, then  we have
\begin{align*}
t^{j-1} \cdot \left( \prod _{i=1}^{k} \beta_{s_i}^{n_i} \right)
&= \beta_{s_{j-1}} \cdots \beta_{s_1}\left( \prod _{i=1}^{k} \beta_{s_i}^{n_i} \right) \beta^{-1}_{s_1} \cdots \beta^{-1}_{s_{j-1}}
\\
&= \left( \prod _{i=j}^{k} \beta_{s_i}^{n_i} \right)\left( \prod_{i=1}^{j-1} \beta^{n_i}_{s_i}\right) \text{,}
\end{align*}
and this product falls under the first case.

Similarly, if $j>1$, $n_i=1 \thickspace \allowbreak \forall i \in \{1, \dots, j\}$ and $n_{j+1}=-1$, then the product
\[
t^{1-j} \cdot \left( \prod _{i=1}^{k} \beta_{s_i}^{n_i} \right)
\]
falls under the first case.
\end{proof}

\section{Definition \ref{qquandle} and Theorem \ref{main1}}
\label{def1}

In this section we verify that Definition \ref{qquandle} really does describe medial quandles, prove Theorem \ref{main1}, and mention some useful properties of the medial quandles $Q(N, (m_i),(X_i))$. 

\begin{proposition}
Suppose $M,N,I, (m_i)$ and $(X_i)$ satisfy the hypotheses of Definition \ref{qquandle}. Then the resulting structure $Q(N, (m_i),(X_i))$ is a medial quandle.
\end{proposition}
\begin{proof}
The first property of Definition \ref{mqdef} is satisfied because every $x \in Q_i$ has $x \triangleright x = m_i - m_i+tx+(1-t)x = x$. 

For the second property of Definition \ref{mqdef}, note that if $x\in Q_i$ and $y \in Q_j$, then $Q_i$ certainly has a well-defined element 
\[
x \triangleright ^{-1} y = t^{-1} \cdot (m_i-m_j+x-(1-t)y)+X_i.
\]
Then
\begin{align*}
\beta_y(x) \triangleright ^{-1} y &= t^{-1} \cdot (m_i-m_j+\beta_y(x)-(1-t)y)+X_i
\\
&= t^{-1} \cdot (m_i-m_j+m_j-m_i +tx+(1-t)y-(1-t)y)+X_i = x 
\end{align*}
and
\begin{align*}
\beta_y(x \triangleright ^{-1} y) &= m_j-m_i+t \cdot (x \triangleright ^{-1} y)+(1-t)y + X_i
\\
&= m_j-m_i+m_i-m_j+x-(1-t)y+(1-t)y +X_i=x.
\end{align*}
These equalities imply that $\beta_y$ is both injective and surjective, and that the inverse function is given by $\beta^{-1}_y (x) = x \triangleright^{-1} y$. 

Now, suppose $w \in Q_h, x\in Q_i,y \in Q_j$ and $z \in Q_k$. Then
\begin{align*}
(w\triangleright x)  \triangleright {} & {}(y \triangleright z)\\
={}&(m_i-m_h+tw+(1-t)x+X_h) \triangleright (m_k-m_j+ty+(1-t)z+X_j)
\\
={} &  m_j-m_h+t(m_i-m_h+tw+(1-t)x) \\
&+ (1-t)(m_k-m_j+ty+(1-t)z)+X_h
\\
={}&(-1-t)m_h+tm_i +tm_j +(1-t)m_k+t^2 w+t(1-t)x \\
&+ t(1-t)y + (1-t)^2z +X_h.
\end{align*}
The third property of Definition \ref{mqdef} is satisfied because the last formula is the same if we interchange $x \leftrightarrow y$ and $i \leftrightarrow j$. \end{proof}

We are now ready to prove Theorem \ref{main1}. Here is a more detailed statement.

\begin{theorem}
Let $Q$ be a medial quandle, with orbits $Q_i$ indexed by a nonempty set $I$. Let $i_0$ be a fixed element of $I$. For each $i \in I$ choose a fixed element $q_i \in Q_i$, let $d_i$ be the elementary displacement $d_i = \beta_{q_i} \beta^{-1}_{q_{i_0}}$, and let $\Fix(q_i) = \{ d \in \Dis(Q) \mid d(q_i)=q_i \}$. Then
\[
Q \cong Q(\Dis(Q),(d_i), (\Fix(q_i)) ).
\]
\end{theorem}
\begin{proof}
Before providing details we remark that notation in $\Dis(Q)$ is unusual for a $\Lambda$-module: addition is composition of displacements, written as multiplication, and scalar multiplication is given by $td =\beta_{q^*}d\beta_{q^*}^{-1}$ for any $q^* \in Q$. We use $q^*=q_{i_0}$ in this proof.

To begin, observe that if $i \in I$ then $\Fix(q_i)$ is a subgroup of $\Dis(Q)$: if $d,d' \in \Fix(q_i)$ then $d^{-1}d'(q_i)=d^{-1}(q_i)=q_i$. Lemma \ref{fix} tells us that $td=d \allowbreak \thickspace \forall d \in \Fix(q_i)$, so it is obvious that $\Fix(q_i)$ is closed under scalar multiplication, and is annihilated by $1-t$.

Let $Q'=Q(\Dis(Q),(d_i), (\Fix(q_i)))$, and for $i \in I$ let $Q'_i$ be the copy of $\Dis(Q)/\Fix(q_i)$ that is the $i$th part of $Q'$.

For each index $i$ there is a natural surjection $F_i:\Dis(Q) \to Q_{q_i}$, defined by $F_i(d)=d(q_i)$. If $d,d' \in \Dis(Q)$ and the cosets $d \Fix(q_i)$ and $d' \Fix(q_i)$ are the same, then $d^{-1}d' \in \Fix(q_i)$, so $F_i(d) = d(q_i) = d(d^{-1}d'(q_i)) =d'(q_i)=F_i(d')$. Therefore $F_i$ induces a well-defined function mapping $Q'_i = \Dis(Q)/\Fix(q_i)$ onto $Q_{q_i}$. Taken together, these induced functions define a surjection $F:Q' \to Q$.

It is easy to see that the surjection $F$ is injective. Suppose $F(x)=F(y) \in Q_{q_i}$. Then $x,y \in Q'_i$, so there are displacements $d,d' \in \Dis(Q)$ with $x=d \Fix(q_i)$ and $y = d' \Fix(q_i)$. Then $d(q_i)=F(x)=F(y)=d'(q_i)$, so $d^{-1}d'(q_i) = d^{-1}d(q_i)=q_i$. It follows that $d^{-1}d' \in \Fix(q_i)$, so the cosets $d \Fix(q_i)$ and $d' \Fix(q_i)$ are the same.

To complete the proof, we show that the bijection $F$ is a quandle map. Suppose $d,d' \in \Dis(Q), x = d \Fix(q_i) \in Q'_i$ and $y = d' \Fix(q_j) \in Q'_j$. Then according to Definition \ref{qquandle}, in the usual additive notation for $\Lambda$-modules we have
\[
F(x \triangleright y)=F(d_j-d_i+tx+(1-t)y+\Fix(q_i) )= F(y-ty+tx+d_j-d_i+\Fix(q_i) ).
\]
Therefore if we switch to multiplicative notation in $\Dis(Q)$, we have
\begin{align*}
F(x & \triangleright y)=F_i(d'-td'+td+d_j-d_i) =F_i(d'(td')^{-1}(td)d_jd^{-1}_i)\\
&= d' \beta_{q_{i_0}} (d')^{-1} \beta_{q_{i_0}}^{-1} \beta_{q_{i_0}} d \beta_{q_{i_0}}^{-1} \beta_{q_j}\beta_{q_{i_0}}^{-1}\beta_{q_{i_0}}\beta_{q_i}^{-1}(q_i)= d' \beta_{q_{i_0}} (d')^{-1} d \beta_{q_{i_0}}^{-1} 
\beta_{q_j}(q_i).
\end{align*}
As $d,d'$ and $\beta_{q_{i_0}}^{-1} 
\beta_{q_j}$ are elements of $\Dis(Q)$, the commutativity of $\Dis(Q)$ implies
\[
F(x \triangleright y) = d' \beta_{q_{i_0}} \beta_{q_{i_0}}^{-1} 
\beta_{q_j} (d')^{-1} d (q_i) = d'
\beta_{q_j} (d')^{-1} d (q_i)= d'((d')^{-1} d (q_i) \triangleright q_j).
\]
As $d'$ is an automorphism of $Q$, it follows that
\[
F(x \triangleright y) = (d'(d')^{-1} d (q_i)) \triangleright (d'(q_j)) = d (q_i) \triangleright d'(q_j) = F(x) \triangleright F(y).
\]
\end{proof}

The quandle $Q(N, (m_i),(X_i))$ is not sensitive to some changes in the choices of the elements $m_i$.

\begin{proposition}
\label{nodiff}
Suppose $M,N,I, (m_i)$ and $(X_i)$ satisfy the hypotheses of Definition \ref{qquandle}. For each $i \in I$ let $n_i$ be any element of $N$, and let $m'_i=m_i+(1-t)n_i$. Then 
\[
Q(N, (m_i),(X_i)) \cong Q(N, (m'_i),(X_i)).
\]
\begin{proof}
Define a function 
\[
f:Q(N, (m_i),(X_i)) \to Q(N, (m'_i),(X_i))
\]
by: if $x \in Q_i$ then $f(x)=x-n_i+X_i \in Q'_i$. 

This function $f$ is certainly bijective. Also, if $x \in Q_i$ and $y \in Q_j$ then 
\begin{align*}
f(x \triangleright y)
= {} & m_j-m_i+tx+(1-t)y-n_i+X_i
\\
= {} & m_j+(1-t)n_j -m_i  -(1-t)n_i+tx-tn_i\\
&+(1-t)y-(1-t)n_j+X_i\\
= {} & m'_j-m'_i+tf(x)+(1-t)f(y)+X_i \\
= {} & f(x) \triangleright'f(y)  \text{,}
\end{align*}
where $\triangleright'$ denotes the operation of $Q(N, (m'_i),(X_i))$. \end{proof}
\end{proposition}

Now, suppose $x,y,z \in Q(N, (m_i),(X_i))$; say $x \in Q_i$, $y \in Q_j$ and $z \in Q_k$. Then the image of $x$ under the elementary displacement $\beta_y \beta_z^{-1}$ is
\begin{align*}
\beta_y \beta_z^{-1}(x)&=(t^{-1} \cdot (m_i-m_k+x-(1-t)z)+X_i)\triangleright y
\\
&=m_j-m_i+tt^{-1} \cdot (m_i-m_k+x-(1-t)z)+(1-t)y+X_i
\\
&=x+m_j-m_k+(1-t)(y-z)+X_i \in Q_i.
\end{align*}
As $m_j-m_k+(1-t)(y-z)$ is independent of $x$, we deduce a very simple description of the displacements of $Q(N, (m_i),(X_i))$.

\begin{proposition}
\label{dispro}
If $d$ is a displacement of $Q=Q(N, (m_i),(X_i))$, then there is an element $n\in N$ such that $d(x)=x+n +X_i \thickspace \forall x \in Q_i$ $\allowbreak \forall i \in I$.
\end{proposition}
\begin{proof}
Arbitrary displacements are compositions of elementary displacements, so the proposition follows from the calculation before the statement. 
\end{proof}

As a small abuse of notation, we use $d(x) \equiv x+n$ to abbreviate $d(x)=x+n +X_i \thickspace \allowbreak \forall x \in Q_i \thickspace \allowbreak \forall i \in I$.
\begin{corollary}
\label{discor}
Let $Q=Q(N, (m_i),(X_i))$, and let $N'$ be the set of $n \in N$ such that the formula $d(x) \equiv x+n$ defines a displacement of $Q$. Then $N'$ is a $\Lambda$-submodule of $N$, and there is a surjective $\Lambda$-linear homomorphism $N' \to \Dis(Q)$ under which the image of $n \in N'$ is the displacement $d_n$ given by $d_n(x) \equiv x+n$. 
\end{corollary}
\begin{proof}
The identity map is certainly a displacement, so $0 \in N'$. Also, $N'$ is closed under subtraction: if $d_1(x) \equiv x+n_1$ and $d_2(x)\equiv x+n_2$, then $d_1d_2^{-1}(x) \equiv x+n_1-n_2$. It follows that $N'$ is a subgroup of $N$, and $n \mapsto d_n$ defines a homomorphism of abelian groups $N' \to \Dis(Q)$. Proposition \ref{dispro} tells us this homomorphism is surjective. 

Now, let $q^*$ be a fixed element of $Q$; say $q^* \in Q_j$. If $d(x) \equiv x+n$,
then for every $i \in I$ and every $x \in Q_i$,
\begin{align*}
(td)(x) &= \beta_{q^*} d \beta_{q^*}^{-1}(x)
\\
&= \beta_{q^*} d(t^{-1} \cdot (m_i-m_j+x-(1-t)q^*)+X_i)
\\
&= \beta_{q^*} (t^{-1} \cdot (m_i-m_j+x-(1-t)q^*)+n+X_i)
\\
&=m_j-m_i+t(t^{-1} \cdot (m_i-m_j+x-(1-t)q^*)+n)+(1-t)q^*+X_i
\\
&=x+tn+X_i.
\end{align*}
Therefore $tn \in N'$, and $td$ is given by $(td)(x) \equiv x+tn$. Similarly, 
\begin{align*}
(t^{-1}d)(x) &= \beta_{q^*}^{-1} d \beta_{q^*}(x)
\\
&= \beta_{q^*}^{-1} d(m_j-m_i+tx+(1-t)q^*+X_i)
\\
&= \beta_{q^*}^{-1} (m_j-m_i+tx+(1-t)q^*+n+X_i)
\\
&=t^{-1} \cdot (m_i-m_j+m_j-m_i+tx+(1-t)q^*+n-(1-t)q^*+X_i)
\\
&=x+t^{-1}n+X_i.
\end{align*}
Therefore $t^{-1}n \in N'$, and $t^{-1}d$ is given by $(t^{-1}d)(x) \equiv x+t^{-1}n$.  \end{proof}

\begin{proposition}
\label{dispro2}
Let $N''$ be the submodule of $N$ generated by \[
\left( \bigcap _{i \in I} X_i \right) \cup(1-t)N \cup \{m_j-m_k \mid j,k \in I \}.
\]
Then $N''=N'$, and $\cap_i X_i$ is the kernel of the map $n \mapsto d_n$ of Corollary \ref{discor}.
\end{proposition}
\begin{proof}
As noted before Proposition \ref{dispro}, if $y \in Q_j$ and $z \in Q_k$ then $\beta_y \beta^{-1}_z(x) \equiv x+m_j-m_k +(1-t)(y-z)$. Taking $y=0+X_j$ and $z=0+X_k$, it follows that $m_j-m_k \in N'$. Taking $j=k$ and $z=0+X_k$, it follows that $(1-t)y \in N'$. Also, if $n \in \cap_i X_i$ then the identity map of $Q$ is a displacement with $d(x) \equiv x+n$. We conclude that $N'' \subseteq N'$.

For the opposite inclusion, recall that every displacement of $Q$ is a composition of elementary displacements. As noted before Proposition \ref{dispro}, every elementary displacement of $Q$ is given by $d(x) \equiv x+n$ for some $n \in N''$. It follows that every displacement of $Q$ is given by $d(x) \equiv x+n_1+ \dots +n_k$ for some $n_1, \dots, n_k \in N''$; then $n= n_1+ \dots +n_k \in N''$, and $d(x) \equiv x+n$. If $d(x) \equiv x+n \equiv x+n'$, where $n \in N''$, then $0+n+X_i=0+n'+X_i \thickspace \allowbreak \forall i \in I$, so $n '-n \in \cap_i X_i$. Therefore $n '-n \in N''$, so $n' \in N''$ too.

Finally, $n \in N'$ is included in the kernel of the map of Corollary \ref{discor} if and only if $d(x) \equiv x+n$ defines the identity map of $Q$. This is the case if and only if $n \in X_i \thickspace \allowbreak \forall i \in I$.
\end{proof}

\section{Notation and Terminology Regarding Links}
\label{linknot}
As mentioned in the introduction, in this paper a \emph{link} is an oriented virtual link. A link is described by its \emph{diagrams}. A link diagram $D$ is derived from a collection $C_1, \dots , C_ \mu$ of closed, oriented, piecewise smooth curves in the plane. These curves may have a finite number of transverse (self-)intersections, called crossings. There are no tangential (self-)intersections. At a classical crossing, two short segments are removed from the underpassing arc, one on each side of the crossing. We use $C(D)$ to denote the set of classical crossings in $D$. A virtual crossing is distinguished by a small circle drawn around it. A link is an equivalence class of link diagrams under the equivalence relation generated by Reidemeister moves involving classical crossings, and detour moves involving virtual crossings. These moves all preserve the significance of $C_1, \dots, C_\mu$, so it is reasonable to speak of a link having \emph{components} $K_1, \dots, K_\mu$ corresponding to the curves $C_1, \dots, C_\mu$. Note that the components $K_1, \dots, K_\mu$ of a link are given with a particular indexing; in general, re-indexing the components will result in a new link not equivalent to the original one.

Removing the short segments at classical crossings cuts each curve $C_i$ into arcs, sometimes called the long arcs of $D$ \cite{MI}. The set of these arcs is denoted $A(D)$. N.b. The arcs incident at a virtual crossing simply pass directly through the crossing; neither arc is cut, and the two arcs are not considered to be joined at the crossing.  

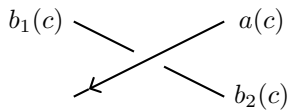
\begin{figure} [bth]
\centering 
\begin{tikzpicture} [>=angle 90]
\draw [thick] [->] (1,.5) -- (-0.8,-0.4);
\draw [thick] (-0.8,-0.4) -- (-1,-.5);
\draw [thick] (-1,.5) -- (-.2,0.1);
\draw [thick] (0.2,-0.1) -- (1,-.5);
\node at (1.5,0.5) {$a(c)$};
\node at (-1.5,0.5) {$b_1(c)$};
\node at (1.5,-0.5) {$b_2(c)$};
\end{tikzpicture}
\caption{The underpassing arcs $b_1(c)$ and $b_2(c)$ of a classical crossing $c$ are on the right and left sides of the overpassing arc $a(c)$, respectively.}
\label{crossfig}
\end{figure}

The \emph{medial quandle} $\MQ(L)$ of a link $L$ is described as follows. If $D$ is a diagram of $L$, then $\MQ(L)$ is the medial quandle generated by elements $q_a$ with $a \in A(D)$, subject to relations stating that at each classical crossing $c$ of $D$ as pictured in Fig.\ \ref{crossfig}, $q_{b_2(c)}=q_{b_1(c)} \triangleright q_{a(c)}$. It is easy to verify that $\MQ(L)$ has one orbit for each component of $L$, with the orbit corresponding to $K_i$ including every $q_a$ such that $a$ is part of the image of $K_i$ in $D$.

The \emph{reduced Alexander module} $\Mr(L)$ of a link $L$ is described as follows. If $D$ is a diagram of $L$, let $\Lambda^{A(D)}$ and $\Lambda^{C(D)}$ be the free $\Lambda$-modules on the sets $A(D)$ and $C(D)$. Define a $\Lambda$-linear map $\varrho_D:\Lambda^{C(D)} \to \Lambda^{A(D)}$ by the formula
\[
\varrho_D(c)=(1-t)a(c)+tb_1(c)-b_2(c) .
\]
Then $\Mr(L) = \coker \varrho_D$. Let $\varsigma_D:\Lambda^{A(D)} \to \Mr(L)$ be the canonical map onto the quotient. For reference, we should mention that the maps $\varrho_D$ and $\varsigma_D$ were denoted $\rho_D \otimes 1$ and $\gamma_D \otimes 1$ in the preceding paper in the series \cite{mvaq4}. We use simpler notation in the present paper because the multivariate version of the Alexander module does not make an appearance here.

Let $A_\mu$ be the $\Lambda$-module obtained from the direct sum $\Lambda \oplus \mathbb Z ^ {\mu-1}$ using the trivial scalar multiplication in the $\mathbb Z$ coordinates. That is, if $\lambda \in \Lambda$ and $n_1, \dots, n_{\mu-1} \in \mathbb Z$ then $t \cdot (\lambda,n_1, \dots, n_{\mu-1}) = (t\lambda,n_1, \dots, n_{\mu-1})$. If $D$ is a diagram of $L=K_1 \cup \dots \cup K_{\mu}$, then there is a function $\kappa_D:A(D) \to \{1, \dots, \mu\}$, with $\kappa_D(a)=i$ if $a$ is part of the image of $K_i$ in $D$. There is also a $\Lambda$-linear map $F:\Lambda^{A(D)} \to A_\mu$ defined by: if $\kappa_D(a)=1$ then $F(a)=(1,0,\dots,0)$; and if $\kappa_D(a)>1$ then $F(a)=(1,0,\dots,0,1,0,\dots,0)$, with the second $1$ in the $\kappa_D(a)$th coordinate. The composition $F \circ \varrho_D$ is identically $0$, so $F$ induces a $\Lambda$-linear map $\Mr(L) \to A_\mu$, which we denote $\phi_\tau$. As explained in \cite{mvaq4}, the notation reflects the fact that $\phi_\tau$ is obtained from a map $\phi_L$ that appears in the multivariate version of Crowell's link module sequence \cite{C1, C3}, by applying a ring homomorphism $\tau$ that simplifies Laurent polynomials in $t_1, \dots, t_\mu$ by setting every $t_i$ equal to $t$. 

Note that the first coordinate of $\phi_\tau$ is the map $\fr$ mentioned in the introduction. It turns out that the kernels of the two maps are closely connected to each other.

\begin{proposition}
\label{oldprop}
The kernel of $\phi_\tau$ is $(1-t) \cdot \ker \fr$.
\end{proposition}
\begin{proof}
The statement of this proposition is the same as the statement of \cite[Prop.\ 9]{mvaq4}. The discussion in \cite{mvaq4} was focused on classical links, but the same proof applies to virtual links too. N.b.\ The maps $\varrho_D$ and $\varsigma_D$ were denoted $\rho_D \otimes 1$ and $\gamma_D \otimes 1$ in \cite{mvaq4}.
\end{proof}

\begin{figure} [bth]
\centering
\begin{tikzpicture} [>=angle 90]
\draw [thick] [->] (1,.5) -- (-0.8,-0.4);
\draw [thick]  (-0.8,-0.4) -- (-1,-.5);
\draw [thick] (-1,.5) -- (-.2,0.1);
\draw [thick] [<-] (0.8,-0.4) -- (0.2,-0.1);
\draw [thick] (0.8,-0.4) -- (1,-.5);
\draw [thick] [->] (5,.5) -- (3.2,-0.4);
\draw [thick]  (3.2,-0.4) -- (3,-.5);
\draw [thick]  [<-] (3.2,.4) -- (3.8,0.1);
\draw [thick]  (3,0.5) -- (3.2,.4);
\draw [thick] (5,-.5) -- (4.2,-0.1);
\node at (0,-1) {$w(c)=+1$};
\node at (4,-1) {$w(c)=-1$};
\end{tikzpicture}
\caption{The writhe of a classical crossing $c$.}
\label{wfig}
\end{figure}
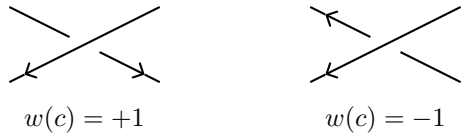

In \cite{peri} we introduced \emph{peripheral elements} in $\Mr(L)$. If $D$ is a diagram of $L=K_1 \cup \dots \cup K_{\mu}$, then for $1 \leq i \leq \mu$ an element $x \in \Mr(L)$ is a \emph{meridian of} $K_i$ in $\Mr(L)$ if $\phi_\tau(x)$ is equal to $\phi_\tau(a)$ for some arc $a \in A(D)$ such that $\kappa_D(a)=i$. The set of meridians of $K_i$ in $\Mr(L)$ is denoted $M_i(L)$. The \emph{longitude} of $K_i$ in $\Mr(L)$ is the element given by the formula 
\[
\chi_i(L) =  \varsigma_D \left (\sum_{\substack{c \in C(D)\\ \kappa_D(b_1(c))=i}} w(c)a(c) - \frac{1}{2} \sum_{\substack{c \in C(D)\\ \kappa_D(b_1(c))=i}} w(c)(b_1(c)+b_2(c)) \right).
\]
Here $w(c)$ denotes the writhe of $c$; see Fig.\ \ref{wfig}. As mentioned in the introduction, these peripheral elements are link invariants, in the sense that if $L$ and $L'$ are representatives of the same link type, then there is an isomorphism $f:\Mr(L) \to \Mr(L')$ with $f(M_i(L))=M_i(L')$ and $f(\chi_i(L))=\chi_i(L')$ for every $i \in \{1, \dots, \mu \}$.

The definition of $\chi_i(L)$ is rather complicated, so it is often convenient to work with a special case. It is easy to see that every link has diagrams with \emph{alternating writhes}, i.e.\ each arc $a \in A(D)$ appears as the underpassing arc of one classical crossing of writhe $-1$ and one classical crossing of writhe $+1$ \cite[Prop.\ 49]{mvaq4}. If $D$ is such a diagram, then the longitudes are given by
\[
\chi_i(L) =  \varsigma_D \left (\sum_{\substack{c \in C(D)\\ \kappa_D(b_1(c))=i}} w(c)a(c) \right).
\]
This simpler formula for the longitudes is our original definition from \cite{mvaq4}. The more complicated general definition was developed so that an arbitrary link diagram can be used to define the longitudes, making it possible to verify invariance under the Reidemeister moves as in \cite{peri}.

\section{Proof of Theorem \ref{main}}
\label{mainproof}

We prove Theorem \ref{main} by defining a quandle map $q_\xi$ from $\MQ(L)$ to the quandle $Q(\ker \fr, m_1, \dots, m_\mu, X_1, \dots, X_\mu)$ mentioned in the theorem, and then showing that $q_\xi$ satisfies the three requirements of Proposition \ref{miso}.

Let $L$ be a link. As mentioned above, $L$ has a diagram $D$ which satisfies the \emph{alternating writhes} requirement, that every $a \in A(D)$ is the underpassing arc at two classical crossings of opposite writhes. 

\emph{Throughout this section we assume $D$ is a diagram with alternating writhes.}
 
We index the arcs and classical crossings of $D$ in the following way. For each $i \in \{1, \dots, \mu \}$, choose an arc $b_{i0} \in A(D)$ that belongs to the link component $K_i$, such that the orientation of $K_i$ has $b_{i0}$ oriented from a classical crossing of writhe $1$ to a classical crossing of writhe $-1$. Start walking along $b_{i0}$, in the direction given by the orientation of $K_i$. When the end of $b_{i0}$ is reached, index that classical crossing as $c_{i0}$, the overpassing arc at that classical crossing as $a_{i0}$, and the next arc of $K_i$ as $b_{i1}$. Continue this process as long as possible, always passing from $b_{ij}$ to $b_{i(j+1)}$ at a classical crossing $c_{ij}$ whose overpassing arc is $a_{ij}$. As $D$ has alternating writhes, $w(c_{ij})$ is always $(-1)^{j+1}$.  Also, $D$ has an even number of arcs belonging to $K_i$; say there are $2n_i$ of them. We consider the second indices of $a_{ij},b_{ij}$ and $c_{ij}$ modulo $2n_i$.
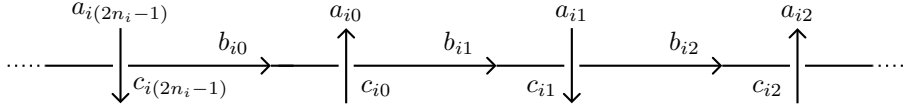
\begin{figure} [t]
\centering
\begin{tikzpicture} [>=angle 90]
\draw [thick] [<-] (-6,-.5) -- (-6,.5);
\draw [thick] [->] (-5.9,0) -- (-4,0);
\draw [thick, dotted] (-7.5,0) -- (-7,0);
\draw [thick] (-7,0) -- (-6.1,0);
\draw [thick] (-4,0) -- (-3.1,0);
\draw [thick] [<-] (-3,.5) -- (-3,-.5);
\draw [thick] (-3.7,0) -- (-3.9,0);
\node at (-5.2,-.3) {$c_{i(2n_i-1)}$};
\node at (-6,.7) {$a_{i(2n_i-1)}$};
\node at (-2.6,-.3) {$c_{i0}$};
\node at (-3,.7) {$a_{i0}$};
\draw [thick] [->] (-2.9,0) -- (-1,0);
\draw [thick] (-1,0) -- (-.1,0);
\draw [thick] [<-] (0,-.5) -- (0,.5);
\node at (-0.4,-.3) {$c_{i1}$};
\node at (0,.7) {$a_{i1}$};
\draw [thick] [->] (.1,0) -- (2,0);
\draw [thick] (2,0) -- (2.9,0);
\draw [thick] [->] (3,-.5) -- (3,.5);
\node at (2.6,-.3) {$c_{i2}$};
\node at (3,.7) {$a_{i2}$};
\draw [thick] (3.1,0) -- (4,0);
\draw [thick, dotted] (4.5,0) -- (4,0);
\node at (-4.5,.3) {$b_{i0}$};
\node at (-1.5,.3) {$b_{i1}$};
\node at (1.5,.3) {$b_{i2}$};
\end{tikzpicture}
\caption{Indexing the arcs and classical crossings of a diagram with alternating writhes.}
\label{indexfig}
\end{figure}

This indexing scheme is illustrated in Fig.\ \ref{indexfig}. Note that if $j$ is even, then $b_1(c_{ij})= b_{i(j+1)}$ and $b_2(c_{ij})=b_{ij}$; while if $j$ is odd, then $b_2(c_{ij})= b_{i(j+1)}$ and $b_1(c_{ij})=b_{ij}$.

For each $i \in \{1, \dots, \mu \}$, define a function $\xi_i:\{b_{i0}, \dots, b_{i(2n_i-1)}\} \to \ker \fr$ as follows. If $0 \leq j<2n_i$ and $j$ is even, then
\[
\xi_i(b_{ij}) = \sum_{k=0}^{j-1} (-1)^{k+1} \varsigma_D(a_{ik}).
\]
In particular, $\xi_i(b_{i0})=0$. On the other hand if $0<j<2n_i$ and $j$ is odd, then
\[
\xi_i(b_{ij}) = t^{-1}\varsigma_D(b_{i0}) +t^{-1} \cdot \sum_{k=0}^{j-1} (-1)^{k+1}\varsigma_D(a_{ik}).
\]

\begin{lemma}
\label{xilem}
If $1 \leq i \leq \mu$ and $0 \leq j \leq 2n_i-1$, then $(1-t)\xi_i(b_{ij})=\varsigma_D(b_{ij}-b_{i0})$.
\end{lemma}
\begin{proof}
Of course if $j=0$ then $(1-t)\xi_i(b_{ij})=0=\varsigma_D(b_{ij}-b_{i0})$.

Suppose $0 < j < 2n_i$, and $j$ is even. Then according to the definition of $\varsigma_D$,
\begin{align*}
(1-t)\xi_i(b_{ij}) =& \sum_{k=0}^{j-1} (-1)^{k+1}(1-t) \varsigma_D(a_{ik})
\\
 =& \sum_{k=0}^{j-1} (-1)^{k+1}\varsigma_D(b_2(c_{ik})-tb_1(c_{ik})).
\end{align*}
As noted above, if $k$ is even then $b_1(c_{ik}) =b_{i(k+1)}$ and $b_2(c_{ik})=b_{ik}$, and if $k$ is odd then $b_2(c_{ik}) =b_{i(k+1)}$ and $b_1(c_{ik})=b_{ik}$. Therefore 
\begin{align*}
(1-t)\xi_i(b_{ij})
=&(-1)^{0+1} \varsigma_D(b_{i0}) + \sum_{\substack{k=2\\ k \text{ even}}}^{j-2}((-1)^{k-1+1}+(-1)^{k+1})\varsigma_D(b_{ik}) 
\\
&+(-1)^{j-1+1} \varsigma_D(b_{ij}) + \sum_{\substack{k=1\\ k \text{ odd}}}^{j-1}((-1)^{k-1+1}+(-1)^{k+1})(-t\varsigma_D(b_{ik}))
\\
=& -\varsigma_D(b_{i0}) +0 +\varsigma_D(b_{ij})+0.
\end{align*}

Now, suppose $0 < j <2n_i$, and $j$ is odd. Then
\begin{align*}
(1-t)\xi_i(b_{ij}) =& (1-t)t^{-1}\varsigma_D(b_{i0}) +t^{-1} \cdot \sum_{k=0}^{j-1} (-1)^{k+1}(1-t)\varsigma_D(a_{ik})
\\
 = {} & (t^{-1}-1)\varsigma_D(b_{i0}) + t^{-1} \cdot \sum_{k=0}^{j-1} (-1)^{k+1}\varsigma_D(b_2(c_{ik})-tb_1(c_{ik})) .
\end{align*}
Recall again that if $k$ is even then $b_1(c_{ik}) =b_{i(k+1)}$ and $b_2(c_{ik})=b_{ik}$, and if $k$ is odd then $b_2(c_{ik}) =b_{i(k+1)}$ and $b_1(c_{ik})=b_{ik}$. Therefore
\begin{align*}
(1-t)\xi_i(b_{ij}) = {} & (t^{-1}-1)\varsigma_D(b_{i0}) +  \sum_{k=0}^{j-1} (-1)^{k+1}\varsigma_D(t^{-1}b_2(c_{ik})-b_1(c_{ik}))
\\
={} &(t^{-1}-1)\varsigma_D(b_{i0}) +(-1)^{0+1}t^{-1}\varsigma_D(b_{i0})\\
&+ \sum_{\substack{k=2\\ k \text{ even}}}^{j-1}t^{-1} \cdot ((-1)^{k-1+1}+(-1)^{k+1})\varsigma_D(b_{ik})
\\
&+\sum_{\substack{k=1\\ k \text{ odd}}}^{j-2}((-1)^{k-1+1}+(-1)^{k+1})(-\varsigma_D(b_{ik}))+(-1)^{j-1+1}(-\varsigma_D(b_{ij}))
\\
=&-\varsigma_D(b_{i0})+0+0+\varsigma_D(b_{ij}).
\end{align*}
\end{proof}

Let $Q_D$ be the quandle $Q(\ker \fr, \varsigma_D(b_{10}), \dots, \varsigma_D(b_{\mu 0}), X_1, \dots, X_\mu)$ given by Definition \ref{qquandle}; here for $1 \leq i \leq \mu$, $X_i$ is the submodule of $\ker \fr$ generated by $\chi_i(L)$. It is convenient to use the specific elements $\varsigma_D(b_{i0})$ in our discussion, but Propositions \ref{nodiff} and \ref{oldprop} tell us that $Q_D$ is not changed (up to quandle isomorphisms) if we replace $\varsigma_D(b_{10}), \dots, \varsigma_D(b_{\mu 0})$ with any other elements $m_1, \dots, m_\mu$ such that $\phi_\tau(m_i)=\phi_\tau(\varsigma_D(b_{i0})) \thickspace \allowbreak \forall i \in \{1, \dots, \mu \}$. Also, the fact that the sets of meridians and the individual longitudinal elements in $\Mr(L)$ are invariant (up to module isomorphisms) under Reidemeister moves \cite{peri} implies that $Q_D$ is actually an invariant of link type (up to quandle isomorphisms).

Recall that $Q_D=Q_1 \cup \dots \cup Q_\mu$, where each $Q_i$ is $\ker \fr / X_i$. The functions $\xi_1, \dots, \xi_ \mu$ define a function $\xi:A(D) \to Q_D$, with $\xi(a)=\xi_i(a)+X_i \in Q_i$ whenever $a \in A(D)$ has $\kappa_D(a)=i$.

\begin{lemma}
\label{xilem1} 
If $1 \leq i \leq \mu$ and $0 \leq j \leq 2n_i-1$, let $d_{ij} = \beta_{\xi(b_{ij})}\beta_{\xi(b_{10})}^{-1} \in \Dis(Q_D)$. Then $d_{ij}$ is given by 
\[
d_{ij}(x) \equiv x+ \varsigma_D(b_{ij}-b_{10}).
\]
That is, $d_{ij} =d_{\varsigma_D(b_{ij}-b_{10})}$ in the notation of Corollary \ref{discor}. 
\end{lemma}
\begin{proof}
Suppose $k \in \{1, \dots, \mu \}$ and $x \in Q_k$. Then
\begin{align*}
d_{ij}(x) = {} &  (x \triangleright^{-1} \xi(b_{10})) \triangleright \xi(b_{ij})
\\
= {} &(t^{-1} \cdot (\varsigma_D(b_{k0})-\varsigma_D(b_{10})+x-(1-t)\xi(b_{10}))+X_k)\triangleright \xi(b_{ij})
\\
= {} & tt^{-1} \cdot (\varsigma_D(b_{k0})-\varsigma_D(b_{10})+x-(1-t)\xi(b_{10}))\\
&+\varsigma_D(b_{i0})-\varsigma_D(b_{k0})+(1-t)\xi(b_{ij})+X_k
\\
= {} & \varsigma_D(b_{k0})-\varsigma_D(b_{10})+x-(1-t)\cdot 0\\
&+\varsigma_D(b_{i0})-\varsigma_D(b_{k0})+(1-t)\xi(b_{ij})+X_k
\\
={} &x+(1-t)\xi(b_{ij})+\varsigma_D(b_{i0})-\varsigma_D(b_{10})+X_k.
\end{align*}
According to Lemma \ref{xilem}, it follows that 
\[
d_{ij}(x)=x+\varsigma_D(b_{ij}-b_{i0})+\varsigma_D(b_{i0})-\varsigma_D(b_{10})+X_k.
\]
\end{proof}

\begin{lemma}
\label{xilem3} 
Suppose $r,s \in \mathbb Z$, $1 \leq i \leq \mu$ and $0 \leq j \leq 2n_i-1$. Let
\[
d_{ij,rs}= (\beta_{\xi(b_{10})}^s d_{ij} \beta_{\xi(b_{10})}^{-s})^r = \beta_{\xi(b_{10})}^s d_{ij}^r \beta_{\xi(b_{10})}^{-s} \in \Dis(Q_D).
\]
Then $ d_{ij,rs}$ is given by 
\[
d_{ij,rs}(x) \equiv x+ rt^s\varsigma_D(b_{ij}-b_{10}).
\]
That is, $d_{ij,rs} =d_{rt^s\varsigma_D(b_{ij}-b_{10})}$ in the notation of Corollary \ref{discor}.
\end{lemma}
\begin{proof}
The proof of Corollary \ref{discor} tells us that if $d=d_n$ is the displacement of $Q_D$ given by $d(x) \equiv x+n$, then for any element $q^* \in Q_D$, $\beta_{q^*}d\beta_{q^*}^{-1}$ is the displacement with $\beta_{q^*}d\beta_{q^*}^{-1}(x) \equiv x+tn$, and $\beta_{q^*}^{-1}d\beta_{q^*}$ is the displacement with $\beta_{q^*}^{-1}d\beta_{q^*}=x+t^{-1}n$. As Lemma \ref{xilem1} tells us that $d_{ij}(x) \equiv x+ \varsigma_D(b_{ij}-b_{10})$, we may apply one of the formulas from the preceding sentence $|s|$ times, with $q^* = \xi(b_{10})$, to conclude that
\[
(\beta_{\xi(b_{10})}^s d_{ij} \beta_{\xi(b_{10})}^{-s})(x) \equiv 
 x+ t^s\varsigma_D(b_{ij}-b_{10}) .
\]
It follows that also
\[
(\beta_{\xi(b_{10})}^s d_{ij} \beta_{\xi(b_{10})}^{-s})^{-1}(x) \equiv x -t^s\varsigma_D(b_{ij}-b_{10}).
\]
Applying one of these last two formuals $|r|$ times, we obtain the formula of the statement.
\end{proof}

\begin{proposition}
\label{xipro1}
The subsets $Q_1, \dots, Q_\mu$ are the orbits of $Q_D$.
\end{proposition}
\begin{proof}
Let $y \in \ker \fr$. Then there are elements $\lambda_{ij} \in \Lambda$, where $1\leq i \leq \mu$ and $0 \leq j \leq 2n_i-1$,  such that 
\[
y = \sum _{i=1}^\mu \sum _{j=0}^{2n_i-1} \lambda_{ij} \varsigma_D(b_{ij}).
\]
As $y \in \ker \fr$, $0=\fr(y)=\sum \lambda_{ij}$. Therefore
\begin{align*}
y &= \sum _{i=1}^\mu \sum _{j=0}^{2n_i-1} \lambda_{ij} \varsigma_D(b_{ij})-\sum _{i=1}^\mu \sum _{j=0}^{2n_i-1} \lambda_{ij} \varsigma_D(b_{10})\\
&= \sum _{i=1}^\mu \sum _{j=0}^{2n_i-1} \lambda_{ij} \varsigma_D(b_{ij}-b_{10}).
\end{align*}

For each choice of indices $i$ and $j$, $\lambda_{ij}$ is a sum $r_1 t^{s_1}+r_2t^{s_2} + \dots + r_pt^{s_p}$ for some integers $r_1, \dots, r_p, s_1, \dots, s_p$. Then according to Lemma \ref{xilem3}, the composition
\[
d_y(ij) = d_{ij,r_1s_1} d_{ij,r_2s_2} \cdots d_{ij, r_ps_p} 
\]
is a displacement with $d_y(ij)(x) \equiv x+\lambda_{ij}\varsigma_D(b_{ij}-b_{10})$. Composing all of these displacements together, we see that 
\[
d_y = \prod_{i=1}^ \mu \prod _{j=0}^{2n_i-1} d_y(ij)
\]
is a displacement of $Q_D$ with $d_y(x) \equiv x+y$.

Now, suppose $1 \leq k \leq \mu$ and $q\in Q_k$. Then $q=y+X_k$ for some $y \in \ker \fr$. According to the formula just given, it follows that $q=d_y(0+X_k)$. We conclude that the orbit of $0+X_k$ includes every element $q \in Q_k$.
\end{proof}

\begin{proposition}
\label{xipro2}
The quandle $Q_D$ is generated by $\xi(A(D))$.
\end{proposition}
\begin{proof}
This proposition follows from the same proof as Proposition \ref{xipro1}, because (a) if $1 \leq k \leq \mu$ then $0+X_k = \xi(b_{k0}) \in \xi(A(D))$  and (b) the displacements $d_{ij,rs}$ that appear in the proof are all compositions of maps $\beta_q^{\pm 1}$, where $q \in \xi(A(D))$.
\end{proof}
\begin{proposition}
\label{xipro}
If $c \in C(D)$ then $\xi(b_2(c)) = \xi(b_1(c)) \triangleright \xi(a(c))$ in $Q_D$.
\end{proposition}
\begin{proof}
Suppose $c=c_{ij} \in C(D)$, and $j$ is even. As noted above, $b_1(c_{ij})= b_{i(j+1)}$ and $b_2(c_{ij})=b_{ij}$. Then we have
\begin{align*}
\xi(b_2(c))&=\xi_i(b_{ij}) +X_i = \sum_{k=0}^{j-1} (-1)^{k+1} \varsigma_D(a_{ik}) +X_i
\\
&= (-1)^j\varsigma_D(a_{ij}) + tt^{-1} \cdot \sum_{k=0}^j(-1)^{k+1} \varsigma_D(a_{ik}) +X_i
\\
&= (-1)^j\varsigma_D(a_{ij})+t \cdot (\xi_i(b_{i(j+1)})- t^{-1} \varsigma_D(b_{i0}))+X_i
\\
&= \varsigma_D(a(c))+t  \xi_i(b_{i(j+1)})-\varsigma_D(b_{i0})+X_i.
\end{align*}
Lemma \ref{xilem} tells us that if $\kappa=\kappa_D(a(c))$, then $\varsigma_D(a(c)) = \varsigma_D(b_{\kappa 0})+(1-t) \xi_{\kappa}(a(c))$. Hence
\[
\xi(b_2(c))=\varsigma_D(b_{\kappa0})-\varsigma_D(b_{i0})+(1-t) \xi(a(c))+t \xi(b_1(c))+X_i.
\]
According to Definition \ref{qquandle}, it follows that $\xi(b_2(c)) = \xi(b_1(c)) \triangleright \xi(a(c))$ in $Q_D$. 

Now, suppose $c=c_{ij} \in C(D)$, $1 \leq j \leq 2n_i-3$, and $j$ is odd. As noted above, $b_1(c_{ij})= b_{ij}$ and $b_2(c_{ij})=b_{i(j+1)}$. Then we have
\begin{align*}
\xi(b_2(c))&=\xi_i(b_{i(j+1)})+X_i = \sum_{k=0}^j (-1)^{k+1} \varsigma_D(a_{ik}) +X_i
\\
&= (-1)^{j+1}\varsigma_D(a_{ij}) + tt^{-1} \cdot \sum_{k=0}^{j-1}(-1)^{k+1} \varsigma_D(a_{ik}) +X_i
\\
&= (-1)^{j+1}\varsigma_D(a_{ij})+t \cdot (\xi_i(b_{ij})- t^{-1} \varsigma_D(b_{i0}))+X_i
\\
&= \varsigma_D(a(c))+t  \xi_i(b_{ij})-\varsigma_D(b_{i0})+X_i.
\end{align*}
Again, it follows from Lemma \ref{xilem} and Definition \ref{qquandle} that if $\kappa=\kappa_D(a(c))$ then
\begin{align*}
\xi(b_2(c))&=\varsigma_D(b_{\kappa 0})+ (1-t) \xi(a(c))+t \xi(b_1(c))-\varsigma_D(b_{i0})+X_i
\\
&=\xi(b_1(c)) \triangleright \xi(a(c)).
\end{align*}

It remains to consider a crossing $c=c_{ij} \in C(D)$ with $j=2n_i-1$. We have $b_1(c)=b_{i(2n_i-1)}$ and $b_2(c)=b_{i0}$, so
\begin{align*}
\xi(b_2(c))&=0+X_i=\chi_i(L)+X_i = \varsigma_D \left (\sum_{\substack{c \in C(D)\\ \kappa_D(b_1(c))=i}} w(c)a(c) \right) + X_i
\\
&=\sum_{k=0}^{2n_i-1}(-1)^{k+1} \varsigma_D(a_{ik}) +X_i=\sum_{k=0}^j(-1)^{k+1} \varsigma_D(a_{ik}) +X_i.
\end{align*}

The argument now proceeds just as it did for odd $j \leq 2n_i-3$. \end{proof}

The next lemma is implicit in \cite{mvaq4}, where the discussion involves the idea of ``augmenting'' $\MQ(L)$ with a group. This idea is not necessary here, so we provide a direct proof for the reader's convenience.

\begin{lemma}
\label{betaed}
There is a $\Lambda$-linear epimorphism $\delta:\ker \fr \to \Dis(\MQ(L))$ such that $\delta \varsigma_D(a-a') = \beta _{q_a} \beta_{q_{a'}}^{-1} \thickspace \allowbreak \forall a,a' \in A(D)$. If $1 \leq i \leq \mu$, then $\delta(\chi_i(L))$ is a displacement whose restriction to the $K_i$ orbit of $\MQ(L)$ is the identity map.
\end{lemma}
\begin{proof}
Choose any fixed element $a^* \in A(D)$. There is certainly a function $\Delta:A(D) \to \Dis(\MQ(L))$ with $\Delta(a)=\beta_{q_a}\beta_{q_{a^*}}^{-1} \thickspace \allowbreak \forall a \in A(D)$. This function defines a $\Lambda$-linear map $\Lambda^{A(D)} \to \Dis(\MQ(L))$, which we also denote $\Delta$. Remember that notation in $\Dis(\MQ(L))$ is unusual for a $\Lambda$-module: addition is composition of displacements, written as multiplication, and scalar multiplication is given by $td =\beta_{q_{a^*}}d\beta_{q_{a^*}}^{-1}$.

Let $c$ be a classical crossing of $D$. Then $q_{b_2(c)}=q_{b_1(c)} \triangleright q_{a(c)}$ in $\MQ(L)$. It follows that if $x \in \MQ(L)$ then 
\begin{align*}
\beta_{q_{b_2(c)}}(x) &= x \triangleright q_{b_2(c)}= x \triangleright (q_{b_1(c)} \triangleright q_{a(c)})
\\
&= ((x \triangleright^{-1} q_{a(c)}) \triangleright q_{a(c)}) \triangleright (q_{b_1(c)} \triangleright q_{a(c)})
\end{align*}
so since $\MQ(L)$ is a medial quandle,
\begin{align*}
\beta_{q_{b_2(c)}}(x) &= ((x \triangleright^{-1} q_{a(c)}) \triangleright q_{b_1(c)}) \triangleright (q_{a(c)} \triangleright q_{a(c)})
\\
&= ((x \triangleright^{-1} q_{a(c)}) \triangleright q_{b_1(c)}) \triangleright q_{a(c)} = \beta_{q_{a(c)}} \beta_{q_{b_1(c)}}\beta_{q_{a(c)}}^{-1}(x).
\end{align*}
That is, $\beta_{q_{b_2(c)}}=\beta_{q_{a(c)}} \beta_{q_{b_1(c)}}\beta_{q_{a(c)}}^{-1}$. Therefore 
\begin{align*}
\Delta(\varrho_D(c))&=\Delta((1-t)a(c)+tb_1(c)-b_2(c))
\\
&=\Delta(a(c)+tb_1(c)-ta(c)-b_2(c))
\\
&=\Delta(a(c))\Delta(tb_1(c))\Delta(-ta(c))\Delta(-b_2(c))
\\
&=\Delta(a(c))(t\Delta(b_1(c)))(t(\Delta(a(c))))^{-1}(\Delta(b_2(c)))^{-1}
\\
&=\beta_{q_{a(c)}}\beta_{q_{a^*}}^{-1}\beta_{q_{a^*}}\beta_{q_{b_1(c)}}\beta_{q_{a^*}}^{-2}(\beta_{q_{a^*}}\beta_{q_{a(c)}}\beta_{q_{a^*}}^{-2})^{-1}(\beta_{q_{b_2(c)}}\beta_{q_{a^*}}^{-1})^{-1}
\\
&=\beta_{q_{a(c)}}\beta_{q_{b_1(c)}}\beta_{q_{a^*}}^{-2}(\beta^2_{q_{a^*}}\beta^{-1}_{q_{a(c)}}\beta_{q_{a^*}}^{-1})(\beta_{q_{a^*}}\beta_{q_{b_2(c)}}^{-1})
\\
&=\beta_{q_{a(c)}}\beta_{q_{b_1(c)}}\beta^{-1}_{q_{a(c)}}\beta_{q_{b_2(c)}}^{-1} =1 \text{,}
\end{align*}
the identity element of $\Dis(\MQ(L))$. It follows that $\Delta$ induces a well-defined $\Lambda$-linear map $\Mr(L) \to \Dis(\MQ(L))$, and this map in turn restricts to a $\Lambda$-linear map $\delta:\ker \fr \to \Dis(\MQ(L))$. If $a,a' \in A(D)$ then $\delta(\varsigma_D(a-a')) = \Delta(a-a') = \Delta(a)(\Delta(a'))^{-1} = \beta_{q_{a}}\beta^{-1}_{q_{a^*}}\beta_{q_{a^*}}\beta_{q_{a'}}^{-1}=\beta_{q_{a}}\beta_{q_{a'}}^{-1}$, as required.

The image of $\delta$ contains $\beta_{q_{a}}\beta_{q_{a'}}^{-1}$ for all choices of $a,a' \in A(D)$. The quandle $\MQ(L)$ is generated by $\{ q_a \mid a \in A(D)\}$, so Lemma \ref{disqlem} tells us that $\delta$ is surjective. 

It remains to prove the second sentence of the statement: the restriction of $\delta (\chi_i(L))$ to the $K_i$ orbit of $\MQ(L)$ is always the identity map. Observe that 
\begin{align*}
\delta(\chi_i(L))(q_{b_{i0}}) &= \delta \left( \sum_{j=0}^{2n_i-1} w(c_{ij})a(c_{ij}) \right)(q_{b_{i0}})
\\
&=\delta(a_{i(2n_i-1)}-a_{i(2n_i-2)}+a_{i(2n_i-3)}- \dots +a_{i1}-a_{i0})(q_{b_{i0}})
\\
&=\beta_{a_{i(2n_i-1)}} \beta^{-1}_{a_{i(2n_i-2)}} \beta_{a_{i(2n_i-3)}} \dots \beta_{a_{i1}}\ \beta_{a_{i0}}^{-1}(q_{b_{i0}})
\\
&=\beta_{a_{i(2n_i-1)}} \beta^{-1}_{a_{i(2n_i-2)}} \beta_{a_{i(2n_i-3)}} \dots \beta_{a_{i1}}(q_{b_{i1}})
\\
&=\beta_{a_{i(2n_i-1)}} \beta^{-1}_{a_{i(2n_i-2)}} \beta_{a_{i(2n_i-3)}} \dots \beta^{-1}_{a_{i2}}(q_{b_{i2}})
\\
&= \dots= \beta_{a_{i(2n_i-1)}} \beta^{-1}_{a_{i(2n_i-2)}}(q_{b_{i(2n_i-2)}})=\beta_{a_{i(2n_i-1)}}(q_{b_{i(2n_i-1)}})=q_{b_{i0}}.
\end{align*}
Thus $\delta(\chi_i(L))$ is a displacement of $\MQ(L)$ that fixes $q_{b_{i0}}$. As $\Dis(\MQ(L))$ is abelian, it follows that for every $d \in \Dis(\MQ(L))$, $\delta(\chi_i(L))$ also fixes $d(q_{b_{i0}})$:
\[
\delta(\chi_i(L))(d(q_{b_{i0}})) = d(\delta(\chi_i(L))(q_{b_{i0}})) = d(q_{b_{i0}}).
\]
The $K_i$ orbit of $\MQ(L)$ is $\{d(q_{b_{i0}}) \mid d \in \Dis(\MQ(L)) \}$, so the proof is complete.
\end{proof}

As mentioned above, Lemma \ref{betaed} is extracted from the discussion in \cite{mvaq4}. The map $\delta$ was denoted $\beta e_D$ there.  We should also mention that in general, $\delta$ is not injective. For example, the discussion of the classical Hopf link $H$ in Sec.\ \ref{exsec} implies that $ \ker\phi_H^{\textup{red}} \cong \Lambda / (1-t)$ is isomorphic to $\mathbb Z$ as an abelian group, while the only element of $\Dis(\MQ(H))$ is the identity map of $\MQ(H)$.

We are now ready to complete the proof of Theorem \ref{main}.

\begin{theorem}
\label{main2}
The function $\xi$ defines an isomorphism $q_\xi:\MQ(L) \to Q_D$, given by $q_\xi(q_a)= \xi(a) \thickspace \allowbreak \forall a \in A(D)$. 
\end{theorem}
\begin{proof}
The existence of the quandle map $q_\xi$ follows immediately from Proposition \ref{xipro} and the fact that $Q_D$ is a medial quandle. The fact that $q_\xi$ is surjective follows from Proposition \ref{xipro2}, so the first property of Proposition \ref{miso} is satisfied.

For the second property of Proposition \ref{miso}, recall that $\MQ(L)$ has one orbit for each link component $K_i$, with the orbit corresponding to $K_i$ including every $q_a$ with $a \in A(D)$ and $\kappa_D(a)=i$. Suppose $x,y \in \MQ(L)$ have different orbits; say the orbits correspond to $K_i$ and $K_j$, where $i\neq j \in \{1, \dots, \mu\}$. Then $q_\xi(x)\in Q_i$ and $q_\xi(y) \in Q_j$. According to Proposition \ref{xipro1}, these are different orbits in $Q_D$.

It remains to verify the third property of Proposition \ref{miso}. Let $\delta:\ker \fr \to \Dis(\MQ(L))$ be the map of Lemma \ref{betaed}. Then for any $a,a'\in A(D)$,
\[
\Dis(q_\xi)(\delta(\varsigma_D(a-a')))= \Dis(q_\xi)(\beta _{q_a} \beta_{q_{a'}}^{-1})=\beta_{q_\xi(q_a)}\beta^{-1}_{q_\xi(q_{a'})}=\beta_{\xi(a)}\beta^{-1}_{\xi(a')}\text{,}
\]
so a calculation just like the proof of Lemma \ref{xilem1} implies that the displacement $\Dis(q_\xi)(\delta(\varsigma_D(a-a')))$ of $Q_D$ is given by
\[
\Dis(q_\xi)(\delta(\varsigma_D(a-a')))(x) \equiv x+\varsigma_D(a-a').
\]
In the notation of Corollary \ref{discor}, $\Dis(q_\xi)(\delta(\varsigma_D(a-a')))=d_{\varsigma_D(a-a')}$. Corollary \ref{discor} tells us that the map $n \mapsto d_n$ is $\Lambda$-linear, so since the elements $\varsigma_D(a-a')$ generate $\ker \fr$, it follows that the composition
\[
\Dis(q_\xi) \delta: \ker \fr \to \Dis (Q_D)
\]
is simply the map $n \mapsto d_n$. 

Now, suppose $x \in \MQ(L)$, $d\in \Dis(\MQ(L))$ and $\Dis(q_\xi)(d)(q_\xi (x))=q_\xi(x)$. As $\delta:\ker \fr \to \Dis(\MQ(L))$ is an epimorphism, there is an $n \in \ker \fr$ such that $d=\delta(n)$. The previous paragraph tells us that $\Dis(q_\xi)(d)=\Dis(q_\xi)(\delta(n))=d_n$. Then $\Dis(q_\xi)(d)(q_\xi (x))=q_\xi(x)$ implies $d_n(q_\xi(x))=n+q_\xi(x) = q_\xi(x)$ in $Q_D$, so if $x$ is an element of the orbit of $\MQ(L)$ corresponding to $K_i$, then $n+q_\xi(x) = q_\xi(x)$ in $Q_i=\ker\fr / X_i$. This requires $n \in X_i$, i.e., $n$ is an integer multiple of $\chi_i(L)$. But then $\delta(n)$ is an integer power of $\delta(\chi_i(L))$, so Lemma \ref{betaed} tells us that the restriction of $d=\delta(n)$ to the $K_i$ orbit of $\MQ(L)$ is the identity map. Therefore $d(x)=x$. \end{proof}

The proof of Theorem \ref{main} is completed by Theorem \ref{main2}. Corollary \ref{knotcor} follows from Theorem \ref{main} and a fact mentioned before Proposition \ref{miso}, that a medial quandle with only one orbit is isomorphic to the standard Alexander quandle on its displacement group.

\section{Semiregular Medial Quandles}
\label{semireg}

If $M$ is a $\Lambda$-module, then there is a medial quandle structure on $M$ given by the operations $x \triangleright y = tx+(1-t)y$ and $x \triangleright^{-1} y = t^{-1}x + (1-t^{-1})y$. We call a medial quandle obtained in this way a \emph{standard Alexander quandle}, or an \emph{affine} quandle. Note that an affine quandle is described by an instance of Definition \ref{qquandle} with $m_i=m_j$ and $X_i=X_j \thickspace \allowbreak \forall i,j \in I$. A quandle is \emph{quasi-affine} if it is isomorphic to a subquandle of an affine quandle. We refer the reader to Jedli\v{c}ka \emph{et al}.\ \cite{JPSZ2} for a thorough account of the properties of quasi-affine quandles.

Suppose $Q=Q(N,(m_i),(X_i))$ is an instance of Definition \ref{qquandle}. Then there is a function $f:Q \to M$ defined as follows: if $x \in Q_i$, then $f(x)=m_i+(1-t)x$. This function is well defined because $(1-t)X_i=0 \thickspace \allowbreak \forall i \in I$. Also, $f$ is a quandle map into the standard Alexander quandle on $M$: if $x \in Q_i$ and $y\in Q_j$ then
\begin{align*}
f(x \triangleright y) &= f(m_j-m_i+tx+(1-t)y+X_i) \\
&= m_i+(1-t)(m_j-m_i+tx+(1-t)y)\\
&=t(m_i+(1-t)x)+(1-t)(m_j+(1-t)y) \\
&= f(x) \triangleright f(y).
\end{align*}
We refer to the image of $f$ as the \emph{natural quasi-affine image} of $Q$, and to $f$ itself as the \emph{natural map}. Note that the natural quasi-affine image of $Q$ is the union $\cup_i (m_i+(1-t)N)$ of cosets of $(1-t)N$ in $M$. If $Q=\MQ(L)$ for a classical link $L$, this union is the quandle denoted $\Qr(L)$ in \cite{mvaq4}, with $N=\ker \fr$. We adopt the same $\Qr(L)$ notation for virtual links.

In a special case, the natural map provides an isomorphism between $\MQ(L)$ and $\Qr(L)$.

\begin{theorem}
\label{samelong}
Suppose $L$ is a link, and every longitude $\chi_i(L)$ generates the same submodule of $\ker \fr$, i.e.\ $X_1 = \dots = X_\mu$. Then $\MQ(L) \cong \Qr(L)$.
\end{theorem}
\begin{proof}For any link $L$, the longitudes generate the module $\ann(1-t) = \{x \in \Mr(L) \mid (1-t)x=0 \}$ \cite[Theorem 9]{peri}. If $L$ has $X_1 = \dots = X_\mu$, then, $X_1 = \dots = X_\mu= \ann(1-t)$. It follows that the natural map $f:\MQ(L) \to \Qr(L)$ is injective on each $Q_i=\ker \fr /X_i$, because $f(x)=m_i+(1-t)x=m_i+(1-t)y=f(y)$ requires $(1-t)(y-x)=0$, and this requires that $x$ and $y$ belong to the same coset modulo $X_i=\ann(1-t)$. Proposition \ref{oldprop} tells us that $\phi_\tau((1-t)x)=0 \thickspace \allowbreak \forall x \in \ker \fr$, so for every $i \in \{1, \dots, \mu \}$, $\phi_\tau(f(x))=\phi_\tau(m_i) \thickspace \allowbreak \forall x \in Q_i$. As $\phi_\tau(m_i) \neq \phi_\tau(m_j)$ when $i \neq j$, it follows that $f$ is injective on all of $\MQ(L)$. \end{proof}

Theorem \ref{samelong} applies vacuously when $\mu=1$; in this special case it is equivalent to Corollary \ref{knotcor}. Theorem \ref{samelong} also applies to every classical 2-component link we have analyzed. They all have the property that the two longitudes are negatives of each other. This property fails for many virtual 2-component links. 

With a little more work we can explain the relationship between $\MQ(L)$ and $\Qr(L)$ in generality.

\begin{definition}
\label{semi}
A quandle $Q$ is \emph{semiregular} if the only displacement of $Q$ with a fixed point is the identity map.
\end{definition}

Semiregularity is a property of all quasi-affine quandles. We recall the easy proof for the reader's convenience.

\begin{proposition}
Quasi-affine quandles are semiregular.
\end{proposition}
\begin{proof}
If $M$ is a $\Lambda$-module with elements $x,y$ and $z$ then 
\[
\beta_z \beta^{-1}_y(x)=t(t^{-1}x + (1-t^{-1})y)+(1-t)z = x+(1-t)(z-y).
\]
Therefore every elementary displacement of a subquandle of the standard Alexander quandle on $M$ is a function of the form $x \mapsto x+m$ for some $m \in M$. A displacement is a composition of elementary displacements, so every displacement of a subquandle of the standard Alexander quandle on $M$ is a function of the form $x \mapsto x+m$ for some $m \in M$. Such a function can have a fixed point only if $m=0$, and in that case the function is the identity map.
\end{proof}

Here is a general description of the relationship between $\MQ(L)$ and $\Qr(L)$.

\begin{theorem}
\label{semimax}
If $L$ is a link, then $\Qr(L)$ is the maximal semiregular image of $\MQ(L)$. That is: if $Q$ is any semiregular medial quandle and $g:\MQ(L) \to Q$ is a surjective quandle map, then $g$ factors through the natural map $f:\MQ(L) \to \Qr(L)$.
\end{theorem}
\begin{proof}
Suppose $Q$ is a semiregular medial quandle, and $g:\MQ(L) \to Q$ is a surjective quandle map. Let $\MQ(L) = Q_1 \cup \dots \cup Q_\mu$ be the decomposition of $\MQ(L) \cong Q(\ker \fr, m_1, \dots, m_\mu, \allowbreak X_1, \dots, X_\mu)$ mentioned in Definition \ref{qquandle}.

Suppose $x\in Q_i$, $y \in Q_j$, and $f(x)=f(y) \in \Qr(L)$. We claim that then $g(x)=g(y) \in Q$. 

To prove the claim, note that according to Proposition \ref{oldprop}, $\phi_\tau((1-t)x)=0=\phi_\tau((1-t)y)$ and hence
\begin{align*}
\phi_\tau(m_i)&=\phi_\tau(m_i+(1-t)x)=\phi_\tau(f(x)) \\
&= \phi_\tau (f(y))=\phi_\tau(m_j+(1-t)y)=\phi_\tau(m_j).
\end{align*}
It follows that $i=j$. Then $f(x)=f(y)$ requires $(1-t)x=(1-t)y$, so $y-x$ is an element of $\ker \fr$ that is annihilated by $1-t$. According to \cite[Theorem 9]{peri}, it follows that
\begin{equation}
\label{diffform}
y-x=\sum_{k=1}^ \mu \lambda_k \chi_k(L)
\end{equation}
for some $\lambda_1, \dots, \lambda_\mu \in \Lambda$. 

Lemma \ref{invlem1} tells us that $\ker \fr$ is generated by $(1-t)\ker \fr \cup \{ m_p-m_q \mid 1 \leq p,q \leq \mu \}$. Proposition \ref{dispro2} then tells us that the submodule of $\ker \fr$ denoted $N'$ in the proposition is the entire module $\ker \fr$. That is, every $n \in \ker \fr$ has the property that there is a displacement $d_n$ of $\MQ(L)$ given by $d_n(z) \equiv z+n$. For $1 \leq k \leq \mu$, let $d_k$ be the displacement of $\MQ(L)$ given by $d_k(z) \equiv z+ \lambda_k\chi_k(L)$. Then (\ref{diffform}) implies that $y=d_1 \dots d_\mu(x)$.

It is obvious that each displacement $d_k$ acts as the identity map on $Q_k$: $Q_k$ is a copy of $\ker \fr / X_k$, and $d_k$ is defined by addition of $\lambda_k \chi_k(L)$, which is an element of $X_k$. Therefore each $d_k$ is a displacement of $\MQ(L)$ that has a fixed point. 

As mentioned a couple of paragraphs before Proposition \ref{miso} in Sec.\ \ref{mq}, the surjective quandle map $g:\MQ(L) \to Q$ defines a surjective homomorphism $\Dis(g):\Dis(\MQ(L)) \to \Dis(Q)$, with the property that $\Dis(g)(d) \circ g = g \circ d$ for every $d \in \Dis(\MQ(L))$. Note that if $z$ is a fixed point of $d$ then $\Dis(g)(d)(g(z)) = g(d(z)) = g(z)$, so $g(z)$ is a fixed point of $\Dis(g)(d)$. It follows that for each $k$, $\Dis(g)(d_k)$ is a displacement of $Q$ that has a fixed point; $Q$ is semiregular, so $\Dis(g)(d_k)$ must be the identity map of $Q$.

The formula $\Dis(g)(d) \circ g = g \circ d$ also applies to the displacement $d=d_1 \dots d_\mu$. As $y=d_1 \dots d_\mu(x)=d(x)$, the formula implies 
\begin{align*}
g(y)&=g(d(x))=\Dis(g)(d)(g(x))
\\
&=\Dis(g)(d_1 \dots d_\mu)(g(x))=\Dis(g)(d_1) \dots \Dis(g)(d_\mu)(g(x)).
\end{align*}
Each $\Dis(g)(d_k)$ is the identity map of $Q$, so we conclude that $g(y)=g(x)$, as claimed.

The claim implies that $g = g' \circ f$ for some function $g':\Qr(L) \to Q$, so the theorem is satisfied. \end{proof}

\section{Involutory Medial Quandles}
\label{invmed}

A quandle is \emph{involutory} if every translation $\beta_y$ is an involution. Every quandle $Q$ has an involutory quotient $Q_{\inv}$; $Q_{\inv}$ is the quandle with a generator $x_{\inv}$ for each $x \in Q$, and defining relations $
x_{\inv} \triangleright y_{\inv} = (x \triangleright y)_{\inv}$ and $ ((x \triangleright y) \triangleright y)_{\inv}= \allowbreak x_{\inv}  \thickspace \allowbreak \forall x,y \in Q$.

\begin{lemma}
\label{invlem}
Let $Q=Q(N, (m_i),(X_i))$ as in Definition \ref{qquandle}. Then 
\[
x_{\inv} = (x+(1-t^2)n)_{\inv} \thickspace \forall x \in Q \thickspace \forall n \in N \text{,}
\]
and
\[
x_{\inv}=(x+p(1+t)(m_j-m_i))_{\inv} \thickspace \forall x \in Q \thickspace \forall i,j \in I \thickspace \forall p\in \mathbb Z.
\]
\end{lemma}
\begin{proof}
Note first that if $x \in Q_k$ then
\begin{align*}
(t^{-2}x)_{\inv}&=(((t^{-2}x) \triangleright(0+X_k)) \triangleright (0+X_k))_{\inv}
\\
&=((t(t^{-2}x)+(1-t)0+X_k)\triangleright (0+X_k))_{\inv}
\\
&=(t \cdot t(t^{-2}x)+(1-t)0+X_k)_{\inv}
=x_{\inv}.
\end{align*}

Now we turn to the first equality of the lemma. Suppose $x \in Q_k$ and $n \in N$. Let $y=n+X_k \in Q_k$. Then 
\begin{align*}
x_{\inv}&=(t^{-2}x)_{\inv}=(((t^{-2}x) \triangleright y) \triangleright y)_{\inv}
\\
&=((t(t^{-2}x)+(1-t)y+X_k)\triangleright y)_{\inv}
\\
&=(t \cdot(t^{-1}x+(1-t)y)+(1-t)y+X_k)_{\inv}
\\
&=(x+(1-t^2)y+X_k)_{\inv}.
\end{align*}
This verifies the first equality.

For the second equality of the lemma, suppose $i,k \in I$ and $x\in Q_k$. Then
\begin{align*}
x_{\inv}&=(t^{-2}x)_{\inv}=(((t^{-2}x)\triangleright(0+X_i))\triangleright(0+X_i))_{\inv}
\\
&=(((m_k-m_i)+t(m_k-m_i+t(t^{-2}x+(1-t)0))+(1-t)0+X_k)))_{\inv}
\\
&=((1+t)(m_k-m_i)+x+X_k)_{\inv}.
\end{align*}
As $x \in Q_k$ is arbitrary we may apply this equality repeatedly, obtaining 
\[
((1+t)(m_k-m_i)+x+X_k)_{\inv}=(2(1+t)(m_k-m_i)+x+X_k)_{\inv}\text{,}
\]
\[
(2(1+t)(m_k-m_i)+x+X_k)_{\inv}=(3(1+t)(m_k-m_i)+x+X_k)_{\inv}\text{,}
\]
and so on. Therefore $x_{\inv}=(p(1+t)(m_k-m_i)+x+X_k)_{\inv}$ for every positive integer $p$. This holds for every $x\in Q_k$, so it also holds if we replace $x$ with $(-p(1+t)(m_k-m_i)+x+X_k)_{\inv}$; hence $(-p(1+t)(m_k-m_i)+x+X_k)_{\inv}=x_{\inv}$ for every positive integer $p$. All in all, we have $x_{\inv}=(p(1+t)(m_k-m_i)+x+X_k)_{\inv}$ for every integer $p$. 

In the previous paragraph $x \in Q_k$ and $i \in I$ are arbitrary, so if $j \in I$ then the equality $x_{\inv}=(p(1+t)(m_k-m_i)+x+X_k)_{\inv}$ holds if we replace $x$ with $x-p(1+t)(m_k-m_j)$, or replace $i$ with $j$. Therefore
\begin{align*}
&(p(1+t)(m_j-m_i)+x+X_k)_{\inv}
\\
&=(p(1+t)(m_k-m_i)-p(1+t)(m_k-m_j)+x+X_k)_{\inv}
\\
&=(-p(1+t)(m_k-m_j)+x+X_k)_{\inv}=x_{\inv}.
\end{align*}
\end{proof}

\begin{proposition}
\label{invpro}
Let $Q=Q(N,(m_i),(X_i))$ as in Definition \ref{qquandle}. Let $S$ be the subgroup of $N$ generated by
\[
(1-t^2)N \cup \{ (1+t)(m_j-m_i) \mid i,j \in I \}. 
\]
If it happens that $S$ is a $\Lambda$-submodule of $N$, then 
\[
Q_{\inv} \cong Q(N/S, (\pi_S( m_i)), (\pi_S (X_i))) \text{,}
\]
where $\pi_S:M \to M/S$ is the canonical map onto the quotient.
\end{proposition}
\begin{proof}
If $S$ is a $\Lambda$-submodule of $M$ then $Q'=Q(N/S, (\pi_S( m_i)), (\pi_S (X_i)))$ is a well-defined instance of Definition \ref{qquandle}. The $\Lambda$-linear map $\pi_S$ defines a surjective quandle map $Q \to Q'$, which we also denote $\pi_S$.

The quandle $Q'$ is involutory, because according to Lemma \ref{invlem}, if $x \in Q'_i$ and $y \in Q'_j$ then 

\begin{align*}
x \triangleright y &- x \triangleright^{-1}y
\\
&=\pi_S(m_j-m_i+tx+(1-t)y-t^{-1} \cdot (m_i-m_j+x-(1-t)y)+X_i)
\\
&=\pi_S(t^{-1}\cdot ((t+1)(m_j-m_i)+(t^2-1)x+(1-t^2)y)+X_i)
\\
&= \pi_S(0+X_i).
\end{align*}

It follows that there is a quandle map $f:Q_{\inv} \to Q'$ with $f(x_{\inv})=\pi_S(x) \thickspace \allowbreak \forall x \in Q$. As $\pi_S$ is surjective, $f$ must also be surjective. 

We claim that $f$ is injective. To verify the claim, suppose $f(x_{\inv})=f(y_{\inv})$. By definition, $f$ maps $(Q_{\inv})_k$ to $Q'_k$ for each index $k$, so $x$ and $y$ belong to the same $Q_k$. Then $\pi_S(x)=\pi_S(y)$ implies that $x$ and $y$ are elements of the same coset of the image of $S$ in $Q_k=N/X_k$. By hypothesis, this image is generated, as an abelian group, by the image of $(1-t^2)N \cup \{ (1+t)(m_j-m_i) \mid i,j \in I \}$. That is, there are integers $p_{ij}$, and an element $n \in N$ such that 
\[
y = x+ (1-t^2)n + \sum _{i,j} p_{ij} (m_j-m_i) + X_k.
\]
Lemma \ref{invlem} now tells us that $x_{\inv} = y_{\inv}.$
\end{proof}

\begin{lemma}
\label{invlem1}
If $L$ is a link then $\ker \fr$ is generated, as an abelian group, by 
\[
(1-t) \ker \fr \cup \{\varsigma_D(b_{i0} - b_{j0}) \mid 1 \leq i,j \leq \mu \}.
\]
It follows that the subgroup of $\ker \fr$ generated by 
\[
(1-t^2) \ker \fr \cup \{ (1+t)\varsigma_D(b_{i0} - b_{j0}) \mid 1 \leq i,j \leq \mu \}
\]
is the $\Lambda$-submodule $(1+t) \ker \fr$.
\end{lemma}
\begin{proof}
Suppose $x \in \ker \fr$, and let $\phi_\tau(x) = (0, p_1, \dots p_{\mu-1}) \in \Lambda \oplus \mathbb Z ^{\mu-1}$. Then 
\[
y= x-\sum_{i=2}^ \mu  p_i \varsigma_D(b_{i0}-b_{10}) \in \ker \phi_\tau.
\]
According to Proposition \ref{oldprop}, it follows that $y \in (1-t) \ker \fr$. The first sentence of the statement follows. For the second sentence of the statement, simply multiply everything in the first sentence by $1+t$. \end{proof}

Proposition \ref{invpro} and Lemma \ref{invlem1} give us the following description of $\IMQ(L)$, the involutory quotient of the medial quandle of $L$.
\begin{theorem}
\label{invthm}
Let $L$ be a link, let $S=(1+t)\ker\fr$, and let $\pi_S:\Mr(L) \to \Mr(L)/S$ be the canonical map onto the quotient. Then $\IMQ(L)$ is isomorphic to the quandle 
\[
Q(\ker \fr/S, \pi_S (m_1), \dots, \pi_S( m_\mu), \allowbreak \pi_S (X_1), \dots, \pi_S( X_\mu))
\]
given by Definition \ref{qquandle}.
\end{theorem}

Taking the quotient $\ker \fr /(1+t)\ker \fr$ has the same effect as setting $t=-1$ in our entire discussion. In particular, the fact that the longitudes of $L$ are elements of $\ker \fr$ annihilated by multiplication by $1-t$ implies that the images of the longitudes in $\ker \fr /(1+t)\ker \fr$ are annihilated by multiplication by $2$. We deduce the following.

\begin{corollary}
\label{invthmcor}
Let $L$ be a link. If $\ker \fr / (1+t) \ker \fr$ is finite, then 
\[
\frac {\mu |\ker \fr / (1+t) \ker \fr| } {2} \leq | \IMQ(L) | \leq \mu |\ker \fr / (1+t) \ker \fr|.
\]
In the special case $\mu=1$, $| \IMQ(L) | = |\ker \fr / (1+t) \ker \fr|$.
\end{corollary}
\begin{proof}
According to Definition \ref{qquandle} and Theorem \ref{invthm}, $\IMQ(L)$ is the union of $\mu$ disjoint sets, each of which is the quotient of $\ker \fr / (1+t) \ker \fr$ by the image of a longitude, which is an element of order $1$ or $2$. If $\mu=1$  then the longitude $\chi_1(L)$ is $0$ \cite{peri}, so $\IMQ(L)$ is simply a copy of $\ker \fr / (1+t) \ker \fr$. \end{proof}

If $L$ is a classical link, then $\ker \fr / (1+t) \ker \fr$ is isomorphic to $H_1(X_2)$, the first homology group of the cyclic double cover $X_2$ of $\mathbb S^3$ branched over $L$. This homology group is of order $|\det(L)|$, the absolute value of the determinant of $L$. (We use the phrase ``of order $0$'' to mean that a group is infinite.)  The images of $\chi_1(L), \dots, \chi_ \mu(L)$ under $\pi_S$ coincide with the longitudes denoted $\lambda_1, \dots \lambda_\mu$ in \cite{mvaq2}.  

The ideas of the present paper allow us to improve some results about involutory medial quandles of classical links proven in \cite{mvaq2}. For instance, it was proven there that if $\mu \geq 2$ and $\det L \neq 0$, then $|\IMQ(L)| \leq \mu | \det L|/2$. This is the opposite of the first inequality of Corollary \ref{invthmcor}, so we deduce an equality. 

\begin{corollary}
\label{invthmcor1}
Let $L$ be a classical link. If $\mu \geq 2$ and $\det L \neq 0$, then 
\[
\frac {\mu |\det L| } {2} = | \IMQ(L) |.
\]
\end{corollary}

In \cite{mvaq2} we also showed that if $L$ is a 2-component classical link, the images of the two longitudes in $\ker \fr / (1+t)\ker \fr$ are equal. It follows that all 2-component classical links satisfy the involutory version of Theorem \ref{samelong}.

\begin{corollary}
\label{invthmcor2}
Let $L$ be a classical link with $\mu \leq 2$. Then $\IMQ(L)$ is isomorphic to the image of $\Qr(L)$ in the quotient module $\Mr(L)/(1+t)\Mr(L)$.
\end{corollary}

This generalizes a result of \cite{mvaq2}, which had the additional hypothesis that $\det L \neq 0$. The image of $\Qr(L)$ in $\Mr(L)/(1+t)\Mr(L)$ is the quandle denoted $Q_A(L)_ \nu$ in \cite{mvaq2}. 

We should emphasize that  Corollaries \ref{invthmcor1} and \ref{invthmcor2} hold only for classical links. In the next section, we show that they both fail for the virtual Hopf link.

\section{Examples}
\label{exsec}
Our first three examples illustrate two features of the theory we have discussed. The first feature was mentioned in the introduction: when $\mu>1$, the enhanced reduced Alexander module $\Me(L)$ is a strictly stronger invariant than the medial quandle $\MQ(L)$. This feature of the theory is easily justified without troubling to work through examples. It is well known that any quandle invariant of a classical link $L$ is unchanged if $L$ is replaced by its orientation-reversed mirror image, and an effect of such a replacement is to multiply all the linking numbers in $L$ by $-1$. As $\Me(L)$ detects the linking numbers \cite[Lemma 8]{peri}, it follows that if any linking number in $L$ is nonzero, $\Me(L)$ will detect the difference between $L$ and its orientation-reversed mirror image. (We should mention that there is an easily fixed error in \cite[Lemma 8]{peri}.)

The second feature of the theory we illustrate here is less obvious: neither $\Me(L)$ nor $\MQ(L)$ contains enough information to determine the corresponding invariants for sublinks of $L$. This feature is illustrated by comparing Whitehead's link with the link $7^2_8$. The same two links were used for the opposite purpose in \cite{mvaq3}, as they exemplify the fact that the fundamental multivariate Alexander quandle $Q_A(L)$ does determine the $Q_A$ quandles of sublinks of $L$.

\subsection{Whitehead's link}

\begin{figure} [bth]
\centering
\begin{tikzpicture} 
\draw [thick] ( -0.9,0.9) -- ( -0.15,0.15);
\draw [thick] ( 0.15,-0.15) -- ( 0.9,-0.9);
\draw [thick] (0.65,0.65) -- (-0.65,-0.65);
\draw [->] [>=angle 90] [thick, domain=-215:-90] plot ({ (1.1)*cos(\x)}, {(1.1)*sin(\x)});
\draw [thick, domain=-90:-55] plot ({ (1.1)*cos(\x)}, {(1.1)*sin(\x)});
\draw [thick, domain=-35:125] plot ({ (1.1)*cos(\x)}, {(1.1)*sin(\x)});
\draw [thick] ( 0.9,0.9) to [out=45, in=0] ( 0,2);
\draw [thick] [->]  [>=angle 90] ( 0.9,-0.9) to [out=-45, in=0] ( 0,-2);
\draw [thick] ( -0.9,0.9) to [out=135, in=180] (0,2);
\draw [thick] (0,-2) to [out=180, in=225] ( -0.9,-0.9);
\node at (1.4,0.3) {$K_2$};
\node at (1.3,1.6) {$K_1$};
\node at (-1.3,-.4) {$a_3$};
\node at (1.3,-0.4) {$a_1$};
\node at (0.6,0.2) {$a_4$};
\node at (-1.2,1.6) {$a_2$};
\node at (-1.2,-1.6) {$a_5$};
\draw [thick] (6+0.9,0.9) -- (6+0.15,0.15);
\draw [thick] (6-0.15,-0.15) -- (6-0.9,-0.9);
\draw [thick] (6-0.65,0.65) -- (6+0.65,-0.65);
\draw  [>=angle 90] [thick, domain=-215:-90] plot ({6-(1.1)*cos(\x)}, {(1.1)*sin(\x)});
\draw [<-] [>=angle 90] [thick, domain=-90:-55] plot ({ 6-(1.1)*cos(\x)}, {(1.1)*sin(\x)});
\draw [thick, domain=-35:125] plot ({6-(1.1)*cos(\x)}, {(1.1)*sin(\x)});
\draw [thick] (6-0.9,0.9) to [out=135, in=180] (6-0,2);
\draw [thick] (6-0.9,-0.9) to [out=225, in=180] (6-0,-2);
\draw [thick] (6+0.9,0.9) to [out=45, in=0] (6-0,2);
\draw [thick] [<-] [>=angle 90] (6-0,-2) to [out=0, in=-45] ( 6+0.9,-0.9);
\node at (7.4,0.3) {$K'_2$};
\node at (7.2,1.6) {$K'_1$};
\node at (4.7,-.4) {$a'_1$};
\node at (7.3,-0.4) {$a'_3$};
\node at (5.4,0.2) {$a'_4$};
\node at (4.8,1.6) {$a'_2$};
\node at (4.8,-1.6) {$a'_5$};
\end{tikzpicture}
\caption{Diagrams $D$ and $D'$ of two versions of Whitehead's link, $W$ and $W'$.}
\label{whfig}
\end{figure}
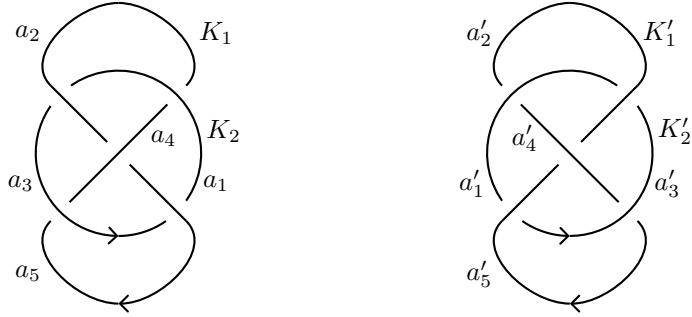

The two versions of Whitehead's link appear in Fig.\ \ref{whfig}. Note that except for the use of apostrophes in $D'$, the diagrams $D$ and $D'$ provide precisely the same functions $a,b_1,b_2$ mapping crossings to arcs, and $\kappa$ mapping arcs to $\{1,2\}$. For instance the arcs in $D$ corresponding to $K_2$ are $a_1$ and $a_3$, and the arcs in $D'$ corresponding to $K'_2$ are $a'_1$ and $a'_3$; an arc that overpasses at a crossing with the same link component is $a_4$ or $a'_4$; at that crossing $b_1$ is $a_5$ or $a'_5$, and $b_2$ is $a_2$ or $a'_2$; a crossing where $a_2$ or $a'_2$ is the overpasser has $b_1$ equal to $a_3$ or $a'_3$ and $b_2$ equal to $a_1$ or $a'_1$; and so on. It follows that for any link invariant determined by the functions $a,b_1,b_2$ and $\kappa$, there is an obvious isomorphism between the invariants of $W$ and $W'$, under which the element corresponding to $a_i$ is matched with the element corresponding to $a'_i$, for each index $i$. Such invariants include Alexander modules, the link module sequences of Crowell \cite{C1, C3}, link groups, and quandles. As we will see, though, the longitudes in $\Me(L)$ are not such invariants.

We proceed to analyze the module $\Mr(W)$. It is generated by the elements $\varsigma_D(a_i)$, with relations corresponding to the crossings of $D$. Three of these crossing relations can be used to describe the other generators in terms of $\varsigma_D(a_3)$ and $\varsigma_D(a_5)$.
\[
\varsigma_D(a_1)=(1-t)\varsigma_D(a_5)+t\varsigma_D(a_3)
\]
\[
\varsigma_D(a_4)=(1-t)\varsigma_D(a_3)+t\varsigma_D(a_5)
\]
\[
\varsigma_D(a_2)=(1-t)\varsigma_D(a_4)+t\varsigma_D(a_5)=(1-t)^2 \varsigma_D(a_3)+(2t-t^2) \varsigma_D(a_5)
\]

The two remaining crossing relations are $\varsigma_D(a_1)=(1-t)\varsigma_D(a_2)+t\varsigma_D(a_3)$ and $\varsigma_D(a_4)=(1-t)\varsigma_D(a_1)+t\varsigma_D(a_2)$. They yield these two relations.
\[
(1-t)\varsigma_D(a_5)+t\varsigma_D(a_3)=(1-t)((1-t)^2 \varsigma_D(a_3)
+(2t-t^2) \varsigma_D(a_5))+t\varsigma_D(a_3)
\]
\begin{align*}
(1-t)\varsigma_D(a_3)+t\varsigma_D(a_5)= {} &(1-t)((1-t)\varsigma_D(a_5)+t\varsigma_D(a_3))
\\
&+t((1-t)^2 \varsigma_D(a_3)+(2t-t^2)\varsigma_D(a_5))
\end{align*}
Both of these relations are equivalent to $(1-t)^3\varsigma_D(a_3) = (1-t)^3\varsigma_D(a_5)$, so 
\[
\Mr(W) \cong \Lambda \oplus (\Lambda/(1-t)^3) \text{,}
\]
with the two direct summands generated by $\varsigma_D(a_5)$ and $\varsigma_D(a_3)-\varsigma_D(a_5)$, respectively. The map $\phi_\tau:\Mr(W) \to \Lambda \oplus \mathbb Z$ is given by
\[
\phi_\tau(x \varsigma_D(a_5) + y(\varsigma_D(a_3)-\varsigma_D(a_5)))=(x,\epsilon(y)),
\]
where $\epsilon:\Lambda \to \mathbb Z$ is the augmentation map given by $\epsilon(t)=1$. 
It follows that 
\[
M_1(W) = \{x \varsigma_D(a_5) + y(\varsigma_D(a_3)-\varsigma_D(a_5)) \mid x=1 \text{ and } \epsilon(y)=0 \}
\]
and
\[
M_2(W) = \{x \varsigma_D(a_5) + y(\varsigma_D(a_3)-\varsigma_D(a_5)) \mid x=1 \text{ and } \epsilon(y)=1 \}.
\]

Of course the same calculations are valid in $\Mr(W')$, with apostrophes. 

Note that $\chi_2(W) = \varsigma_D(a_5-a_2)=(1-t)^2 \varsigma_D(a_5-a_3)$ and $\chi_2(W')=\varsigma_{D'}(a'_2-a'_5)=(1-t)^2 \varsigma_{D'}(a'_3-a'_5)$. We claim that there is no isomorphism $f:\Mr(W) \to \Mr(W')$ with $f(M_1(W))=M_1(W')$, $f(M_2(W))=M_2(W')$ and $f(\chi_2(W))=\chi_2(W')$. The reason is simple: if $f(M_1(W))=M_1(W')$ and $f(M_2(W))=M_2(W')$, then $f(\chi_2(W))=(1-t)^2 f(\varsigma_D(a_5-a_3))$ where $f(\varsigma_D(a_5-a_3))$ is an element of $M_1(W')$ minus an element of $M_2(W')$, while $\chi_2(W') = (1-t)^2 \varsigma_D(a'_3-a'_5)$ where $\varsigma_D(a'_3-a'_5)$ is an element of $M_2(W')$ minus an element of $M_1(W')$. We deduce the property mentioned earlier, that $W$ and $W'$ are distinguished by their enhanced reduced Alexander modules.

For the purpose of comparison with our next example, note also that 
\begin{align*}
\chi_1(W)=&-\varsigma_D(a_1)-\varsigma_D(a_4)+\varsigma_D(a_3)+\varsigma_D(a_2)
\\
=&-(1-t)\varsigma_D(a_5)-t\varsigma_D(a_3)-(1-t)\varsigma_D(a_3)-t\varsigma_D(a_5)
\\
&+\varsigma_D(a_3)+(1-t)^2 \varsigma_D(a_3)+(2t-t^2) \varsigma_D(a_5)
\\
= {} &(1-t)^2\varsigma_D(a_3-a_5)=-\chi_2(W)
\end{align*}
and similarly, 
\begin{align*}
\chi_1(W')&=\varsigma_{D'}(a'_1)+\varsigma_{D'}(a'_4)-\varsigma_{D'}(a'_3)-\varsigma_{D'}(a'_2)
\\
&=(1-t)^2\varsigma_{D'}(a'_5-a'_3)=-\chi_2(W').
\end{align*}

\subsection{The link $7^2_8$}
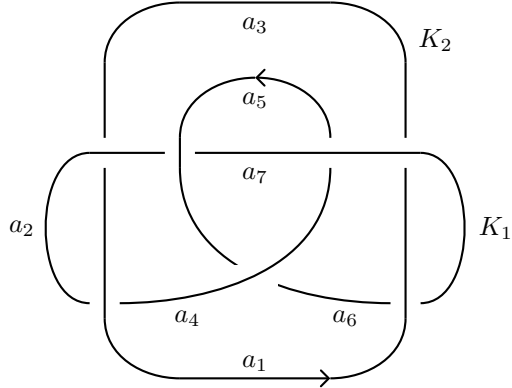
\begin{figure} [bth]
\centering
\begin{tikzpicture} 
\draw [thick] (-1,6) -- (1,6);
\draw [thick] (1,6) to [out=0, in=90] (2,5.2);
\draw [thick] (2,5.2) -- (2,4.2);
\draw [thick] (-1,6) to [out=180, in=90] (-2,5.2);
\draw [thick] (-2,5.2) -- (-2,4.2);
\draw [thick] (-0.8,4) -- (2.2,4);
\draw [thick] (2,3.8) -- (2,1.8);
\draw [thick] (2.2,4) to [out=0, in=0] (2.2,2);
\draw [thick] (1,1) to [out=0, in=-90] (2,1.8);
\draw [thick] (-1,1) to [out=180, in=-90] (-2,1.8);
\draw [thick] [->] [>=angle 90] (-1,1) -- (1,1);
\draw [thick] (-2,3.8) -- (-2,1.8);
\draw [thick] [->] [>=angle 90] (1,4.2) to [out=90, in=0] (0,5);
\draw [thick] (-1,4.2) to [out=90, in=180] (0,5);
\draw [thick] (-1,4.2) -- (-1,3.8);
\draw [thick] (-1,3.8) to [out=-90, in = 180] (1.8,2);
\fill [white] (0.3,2.5) rectangle (-0.6,1.5);
\draw [thick] (1,3.8) to [out=-90, in = 0] (-1.8,2);
\draw [thick] (-1.2,4) -- (-2.2,4);
\draw [thick] (-2.2,4) to [out=180, in=180] (-2.2,2);
\node at (2.4,5.5) {$K_2$};
\node at (3.2,3) {$K_1$};
\node at (0,1.2) {$a_1$};
\node at (-3.1,3) {$a_2$};
\node at (0,5.7) {$a_3$};
\node at (-0.9,1.8) {$a_4$};
\node at (0,4.7) {$a_5$};
\node at (1.2,1.8) {$a_6$};
\node at (0,3.7) {$a_7$};
\end{tikzpicture}
\caption{A diagram $E$ of the link $L=7^2_8$.}
\label{otherfig}
\end{figure}

Our next example is the link $7^2_8$, pictured in Fig.\ \ref {otherfig}. It turns out that $\Mr(7^2_8)$ is generated by the two elements $\varsigma_E(a_1)$ and $\varsigma_E(a_7)$. To verify this, note first that crossing relations give $\varsigma_E(a_1)=(1-t)\varsigma_E(a_7)+t\varsigma_E(a_3) = (1-t)\varsigma_E(a_2)+t\varsigma_E(a_3)$, and hence $(1-t)\varsigma_E(a_7) = (1-t)\varsigma_E(a_2)$. We now have the following equalities.
\begin{align*}
\varsigma_E(a_6) &=(1-t)\varsigma_E(a_1)+t\varsigma_E(a_7) \\
\varsigma_E(a_5)& =(1-t)\varsigma_E(a_4)+t\varsigma_E(a_6) \\
& = (1-t)((1-t)\varsigma_E(a_1)+t\varsigma_E(a_2))+t\varsigma_E(a_6) \\
&=(1-t)^2\varsigma_E(a_1)+t(1-t)\varsigma_E(a_2)+t((1-t)\varsigma_E(a_1)+t\varsigma_E(a_7))\\
&=(1-t)^2\varsigma_E(a_1)+t(1-t)\varsigma_E(a_7)+t((1-t)\varsigma_E(a_1)+t\varsigma_E(a_7))\\
&=(1-t)\varsigma_E(a_1)+t\varsigma_E(a_7) = \varsigma_E(a_6)\\
\varsigma_E(a_4)&=(1-t)\varsigma_E(a_7)+t\varsigma_E(a_5)\\
&=(1-t)\varsigma_E(a_7)+t((1-t)\varsigma_E(a_1)+t\varsigma_E(a_7))\\
&=(t-t^2)\varsigma_E(a_1)+(1-t+t^2)\varsigma_E(a_7)\\
\varsigma_E(a_3)&=(1-t^{-1})\varsigma_E(a_7)+t^{-1}\varsigma_E(a_1)
\\
\varsigma_E(a_2)&=(1-t^{-1})\varsigma_E(a_1)+t^{-1}\varsigma_E(a_4)\\
&=(2-t-t^{-1})\varsigma_E(a_1)+(t^{-1}-1+t)\varsigma_E(a_7)
\end{align*}

Two crossing relations remain. One is
\begin{align*}
\varsigma_E(a_1)= {} & (1-t)\varsigma_E(a_2)+t\varsigma_E(a_3)\\
    = {} & (1-t)((2-t-t^{-1})\varsigma_E(a_1)+(t^{-1}-1+t)\varsigma_E(a_7))\\
    &+t((1-t^{-1})\varsigma_E(a_7)+t^{-1}\varsigma_E(a_1))\\
    = {} & (4-t^{-1}-3t+t^2)\varsigma_E(a_1)+(t^{-1}-3+3t-t^2)\varsigma_E(a_7).
\end{align*}
The other is $\varsigma_E(a_5)=(1-t)\varsigma_E(a_4)+t\varsigma_E(a_6)$, and since $\varsigma_E(a_5)=\varsigma_E(a_6)$, this relation is equivalent to $(1-t)\varsigma_E(a_6)=(1-t)\varsigma_E(a_4)$, or
\[
(1-t)((1-t)\varsigma_E(a_1)+t\varsigma_E(a_7)) 
\]
\[
= (1-t)((t-t^2)\varsigma_E(a_1)+(1-t+t^2)\varsigma_E(a_7)).
\]
This in turn is equivalent to 
\[
(1-t)^2(1-t)\varsigma_E(a_1)=(1-t)(-t+1-t+t^2)\varsigma_E(a_7).
\]
We see that both of these remaining crossing relations are equivalent to $(1-t)^3\varsigma_E(a_1) = (1-t)^3\varsigma_E(a_7)$, so we conclude that 
\[
\Mr(7^2_8) \cong \Lambda \oplus (\Lambda/(1-t)^3) \text{,}
\]
with the two summands generated by $\varsigma_E(a_7)$ and $\varsigma_E(a_1)-\varsigma_E(a_7)$, respectively. The map $\phi_\tau:\Mr(7^2_8) \to \Lambda \oplus \mathbb Z$ is given by
\[
\phi_\tau(x \varsigma_D(a_7) + y(\varsigma_D(a_1)-\varsigma_D(a_7))=(x,\epsilon(y)) .
\]
It follows that 
\[
M_1(7^2_8) = \{x \varsigma_D(a_7) + y(\varsigma_D(a_1)-\varsigma_D(a_7)) \mid x=1 \text{ and } \epsilon(y)=0 \}
\]
and
\[
M_2(7^2_8) = \{x \varsigma_D(a_7) + y(\varsigma_D(a_1)-\varsigma_D(a_7)) \mid x=1 \text{ and } \epsilon(y)=1 \}.
\]

The first longitude of $7^2_8$ is
\begin{align*}
\chi_1(7^2_8) = {} & \varsigma_E(a_1)-\varsigma_E(a_7)+\varsigma_E(a_5)-\varsigma_E(a_4)+\varsigma_E(a_6)-\varsigma_E(a_1)\\
&+\varsigma_E(a_7)-\varsigma_E(a_5)\\
=&-\varsigma_E(a_4)+\varsigma_E(a_6)
\\
= {} &-(t-t^2)\varsigma_E(a_1)-(1-t+t^2)\varsigma_E(a_7)+(1-t)\varsigma_E(a_1)+t\varsigma_E(a_7)
\\
= {} & (1-2t+t^2)\varsigma_E(a_1)-(1-2t+t^2)\varsigma_E(a_7).
\end{align*}

The second longitude of $7^2_8$ is 
\[
\chi_2(7^2_8)=\varsigma_E(a_2)-\varsigma_E(a_7) = (2-t-t^{-1})\varsigma_E(a_1)+(t^{-1}-2+t)\varsigma_E(a_7)
\]
\[
=t^{-1} (1-t)^2(\varsigma_E(a_7)-\varsigma_E(a_1)).
\]
As $(1-t)\chi_2(7^2_8)=0$, it follows that  $\chi_2(7^2_8)=t\chi_2(7^2_8)=(1-t)^2(\varsigma_E(a_7)-\varsigma_E(a_1))=-\chi_1(7^2_8)$.

These calculations show that there is an isomorphism $f:\Mr(W) \to \Mr(7^2_8)$ given by $f(\varsigma_{D}(a_5))=\varsigma_E(a_7)$ and $f(\varsigma_{D}(a_3))=\varsigma_E(a_1)$, and this isomorphism has $f(M_1(W))=M_1(7^2_8)$, $f(M_2(W))=M_2(7^2_8)$, $f(\chi_1(W)) = \chi_1(7^2_8)$ and $f(\chi_2(W)) = \chi_2(7^2_8)$.

We see that the enhanced reduced Alexander module does not distinguish $W$ from $7^2_8$. Theorem \ref{main} tells us that the medial quandle also  does not distinguish $W$ from $7^2_8$. Note that the components of $W$ are both unknots, while one component of $7^2_8$ is a trefoil. The trefoil's involutory medial quandle has cardinality 3, and the unknot's involutory medial quandle has cardinality 1. Therefore the involutory medial quandle distinguishes the trefoil from the unknot; it follows that the medial quandle also distinguishes these two knots. Of course the reduced Alexander module distinguishes them too, as the trefoil and unknot have different Alexander polynomials. We conclude that in general, $\Me(L)$ does not determine the $\Me$ modules of sublinks of $L$, and $\MQ(L)$ does not determine the $\MQ$ quandles of sublinks of $L$. 

This conclusion contrasts with the main result of \cite{mvaq3}, that a quandle $Q_A(L)$ contained in the multivariate Alexander module of $L$ does determine the $Q_A$ quandles of sublinks of $L$. Of course $Q_A(L)$, as a quandle invariant, cannot distinguish $W$ from $W'$. Taken all together, then, these examples indicate that if it is possible to extend the peripheral machinery of \cite{peri} to the multivariate Alexander module in such a way as to incorporate $Q_A(L)$, the resulting invariants might well be stronger than $\Me(L)$, $\MQ(L)$ and $Q_A(L)$. 

\subsection{The classical and virtual Hopf links}

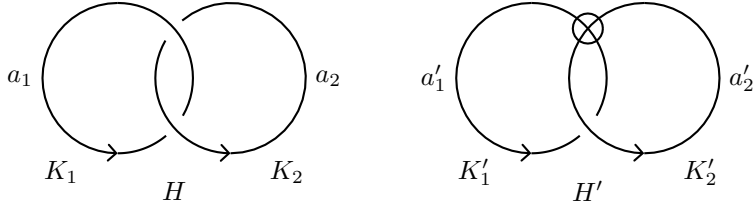
\begin{figure} 
\centering
\begin{tikzpicture} [>=angle 90]
\draw [->] [thick, domain=-270:-90] plot ({-3.5+cos(\x)}, {sin(\x)});
\draw [thick, domain=-50:-90] plot ({-3.5+cos(\x)}, {sin(\x)});
\draw [thick, domain=-30:90] plot ({-3.5+cos(\x)}, {sin(\x)});
\draw [->] [thick, domain=-210:-90] plot ({-2+cos(\x)}, {sin(\x)});
\draw  [thick, domain=-90:130] plot ({-2+cos(\x)}, {sin(\x)});
\draw [->] [thick, domain=-270:-90] plot ({3.5+cos(\x)}, {sin(\x)});
\draw [thick, domain=-90:90] plot ({3.5+cos(\x)}, {sin(\x)});
\draw [->] [thick, domain=-270:-90] plot ({2+cos(\x)}, {sin(\x)});
\draw  [thick, domain=-90:-50] plot ({2+cos(\x)}, {sin(\x)});
\draw  [thick, domain=-30:90] plot ({2+cos(\x)}, {sin(\x)});
\draw  [thick, domain=0:360] plot ({2.75+.2*cos(\x)}, {.66+.2*sin(\x)});
\node at (-4.8,0) {$a_1$};
\node at (-0.7,0) {$a_2$};
\node at (0.7,0) {$a'_1$};
\node at (4.8,0) {$a'_2$};
\node at (2.75,-1.5) {$H'$};
\node at (1.25,-1.25) {$K_1'$};
\node at (4.25,-1.25) {$K_2'$};
\node at (-4.25,-1.25) {$K_1$};
\node at (-1.25,-1.25) {$K_2$};
\node at (-2.75,-1.5) {$H$};
\end{tikzpicture}
\caption{The classical and virtual Hopf links.}
\label{hopfs}
\end{figure}

We finish with two very simple examples to illustrate some differences between classical and virtual links. For convenience, we use $H$ and $H'$ to denote the diagrams pictured in Fig.\ \ref{hopfs} and also the links they represent.

Crossing relations in $H$ imply that $\Mr(H) \cong \Lambda \oplus (\Lambda/(1-t))$, with the two direct summands generated by $\varsigma_H(a_1)$ and $\varsigma_H(a_2-a_1)$, respectively. The module $\Mr(H')$ is described in exactly the same way, with apostrophes attached. The obvious isomorphism between $\Mr(H)$ and $\Mr(H')$ is compatible with the two $\phi_\tau$ maps, so it maps $\Qr(H)$ to $\Qr(H')$. It follows that these two quandles are isomorphic. Similar reasoning leads to the conclusion that $Q_A(H)_\nu$ and $Q_A(H')_\nu$ are isomorphic too. Theorem \ref{samelong} and Corollary \ref{invthmcor2} might lead one to guess that the medial and involutory medial quandles of $H$ are also isomorphic to those of $H'$. In fact, however, $\MQ(H')$ and $\IMQ(H')$ are both strictly larger than $\MQ(H)$ and $\IMQ(H)$.

For $H$, the crossing relations $q_{a_1} \triangleright q_{a_2}=q_{a_1}$ and $q_{a_2} \triangleright q_{a_1} = q_{a_2}$ imply that the medial quandle has only two elements, $q_{a_1}$ and $q_{a_2}$, and $\beta_{q_{a_1}}$ and $\beta_{q_{a_2}}$ are both the identity map. The involutory medial quandle of $H$ is the same.

For $H'$, though, we have only one crossing relation, $q_{a'_1} \triangleright q_{a'_2} = q_{a'_1}$. Recalling the formula $\beta_{ w \triangleright z}=\beta_z \beta_w \beta_z^{-1}$ that appears in the proof of Lemma \ref{disqlem}, we deduce that $\beta_{q_{a'_1}}=\beta_{q_{a'_2}} \beta_{q_{a'_1}} \beta_{q_{a'_2}}^{-1}$. That is, $\beta_{q_{a'_1}}$ and $\beta_{q_{a'_2}}$ commute with each other. It follows that the orbit of $q_{a'_1}$ has only one element, in either $\MQ(H')$ or $\IMQ(H')$. There is no constraint on $\beta_{q_{a'_1}}(q_{a'_2})$, however, so the orbit of $q_{a'_2}$ in $\MQ(H')$ contains infinitely many distinct elements $\beta^n_{q_{a'_1}}(q_{a'_2})$ with $n \in \mathbb Z$, and the orbit of $q_{a'_2}$ in $\IMQ(H')$ contains two distinct elements, $q_{a'_2}$ and $q_{a'_2} \triangleright q_{a'_1}$.

\end{document}